\documentclass[11pt]{article}
\usepackage[a4paper]{anysize}\marginsize{2.5cm}{2.5cm}{1.3cm}{3cm}
\pdfpagewidth=\paperwidth \pdfpageheight=\paperheight
\usepackage{amssymb,amsmath,graphicx}
\usepackage{tikz}
\usepackage{multicol}
\usetikzlibrary{arrows}
\usepackage[latin1]{inputenc}

\let\origLarge=\Large
\renewcommand\LARGE{\bfseries\origLarge}
\makeatletter
\let\orig@item=\@item \def\@item[#1]{\orig@item[\rm #1]}
\renewenvironment{abstract}{\begin{quote}\footnotesize\textbf{\abstractname.}}
{\end{quote}\bigskip}
\newcommand\grant[1]{{\renewcommand\thefootnote{}\footnotetext{#1.}}}
\renewcommand\@seccntformat[1]{\csname the#1\endcsname.\enspace}
\renewcommand\section{\@startsection{section}{1}{\z@}
{-2\baselineskip plus 0.5\baselineskip minus 0.5\baselineskip}
{0.75\baselineskip plus 0.5\baselineskip}{\normalsize\bfseries\centering}}
\renewcommand\subsection{\@startsection{subsection}{2}{\z@}
{-2\baselineskip plus 0.5\baselineskip minus 0.5\baselineskip}
{0.75\baselineskip plus 0.5\baselineskip}{\normalsize\bfseries\centering}}
\renewcommand\paragraph{\@startsection{paragraph}{4}{\z@}{1\baselineskip}
{-0.5em}{\normalsize\bfseries}}
\let\origcaption=\caption \renewcommand\caption[1]{\parbox{0.66\textwidth}
{\origcaption{#1}}}

\renewcommand\@begintheorem[2]{\trivlist\item[\hskip\labelsep
{\bfseries#1 #2.}]\it}
\renewcommand\@opargbegintheorem[3]{\trivlist\item[\hskip\labelsep
{\bfseries#1 #2}] {\bfseries(#3).}\enspace\it\ignorespaces}
\makeatother

\newtheorem{satz}{Satz}[section]
\makeatletter\@addtoreset{equation}{satz}\makeatother

\newtheorem{theorem}[satz]{Theorem}

\newtheorem{lemma}[satz]{Lemma}
\newtheorem{proposition}[satz]{Proposition}
\newtheorem{corollary}[satz]{Corollary}

\newenvironment{definition}[1][\addtocounter{satz}{1}\bf Definition \thesatz]{\trivlist\item[\hskip\labelsep{\it #1.}]}{}
\newenvironment{remark}[1][\addtocounter{satz}{1}\bf Remark \thesatz]{\trivlist\item[\hskip\labelsep{\it #1.}]}{}
\newenvironment{example}[1][\addtocounter{satz}{1}\bf Example \thesatz]{\trivlist\item[\hskip\labelsep{\it #1.}]}{}

\newenvironment{proof}[1][Proof]{\trivlist\item[\hskip\labelsep{\it #1.}]}
{\hspace*{\fill}$\Box$\endtrivlist}
  
\newtheorem{algorithm}[satz]{Algorithm}
\newenvironment{lines}
   {\newcommand\+{\hspace*{2em}\ignorespaces} 
   \begin{list}{}{\leftmargin=\leftmargini\parsep=0cm}\item\begin{obeylines}}
   {\end{obeylines}\end{list}}

\newcommand\subjclass[1]{{\renewcommand\thefootnote{}\footnotetext{2000 
\textit{Mathematics Subject Classification:} #1.}}}
\newcommand\keywords[1]{{\renewcommand\thefootnote{}\footnotetext
{\textit{Keywords.} #1.}}}

\newcommand\engqq[1]{``#1''}

\renewcommand\emptyset{\varnothing}  
\renewcommand\ge{\geqslant}  
\renewcommand\le{\leqslant}  
\renewcommand\geq{\geqslant}  
\renewcommand\leq{\leqslant}  
\renewcommand\epsilon{\varepsilon}
\renewcommand\phi{\varphi}
\renewcommand\bar{\overline}

\renewcommand\tilde{\widetilde}
\renewcommand\O{{\cal O}}
\renewcommand\P{\mathbb P}

\newcommand\be{\begingroup\arraycolsep=0.13888em\begin{eqnarray*}}
\newcommand\ee{\end{eqnarray*}\endgroup}
\renewcommand\to{\longrightarrow}

\renewcommand\mapsto{\mapstochar\longrightarrow}

\newcommand\set[1]{\left\{#1\right\}}

\newcommand\with{\ \vrule\ }

\newlength\matrcolsep \matrcolsep=\arraycolsep

\newcommand\tline{\noalign{\vskip0.4ex}\hline\noalign{\vskip0.65ex}}

\newcommand\converges[2]{\mathop{\longrightarrow}\limits_{n\to\infty}}

\newcommand\N{\mathbb N}
\newcommand\Q{\mathbb Q}
\newcommand\R{\mathbb R}

\newcommand\Z{\mathbb Z}

\newcommand\newop[2]{\newcommand#1{\mathop{\rm #2}\nolimits}}
\newop\mult{mult}
\newop\NS{NS}
\newop\Amp{Amp}
\newop\Pic{Pic}
\newop\Bl{Bl} 
\newop\End{End} 
\newop\Nef{Nef}
\newop\Mov{Mov}
\newop\vol{vol}
\newop\codim{codim}
\newop\BigCone{Big} 
\newop\interior{int}

\newop{\rang}{rang}
\newop{\dv}{div}
\newop{\id}{id}
\newop{\ord}{ord}
\newop{\Div}{Div}
\newop{\Proj}{Proj}
\newop{\im}{im}
\newop{\Face}{Face}
\newop{\conv}{conv}
\newop{\exc}{exc}
\newop{\dom}{dom}
\newop{\Supp}{Supp}
\newop{\Vol}{Vol}
\newop{\Neg}{Neg}
\newop{\Cox}{Cox}

\newop{\Null}{Null}
\newop{\NE}{\overline{{NE}}}
\newop{\Eff}{\overline{{Eff}}}

\begin{document}

   \title{Minkowski decomposition and generators of the moving cone
	for toric varieties}
   \author{Piotr Pokora, David Schmitz, Stefano Urbinati}
   \date{}
   \maketitle
   \thispagestyle{empty}
   \subjclass{Primary 14C20; Secondary 14M25}
   \keywords{Okounkov body, Toric variety, Movable cone, Secondary fan}
   
\grant{The first author was partially supported by WCMCS, Warsaw}
\grant{The second author was supported by DFG grant BA 1559/6-1}

\thispagestyle{empty}
\begin{abstract}
We prove that for smooth projective toric varieties the Okounkov body of a $T$-invariant pseudo-effective divisor with respect to a $T$-invariant flag decomposes as a finite Minkowski sum of indecomposable polytopes and that the set of these polytopes correspond to a finite Minkowski basis whose elements span the extremal rays in the secondary fan. In fact, 
the Minkowski basis does not depend on the choice of the $T$-invariant flag. 
Moreover, we present an algorithm which computes the Minkowski basis.

\end{abstract}


\section*{Introduction}
The definition of Okounkov bodies originates in papers due to A. Okounkov from the middle of the 1990s in the context of representation theory (for instance \cite{ok}). More recently, Lazarsfeld and Musta{\c{t}}{\u{a}} \cite{lm} and independently Kaveh and Khovanskii \cite{kk} initiated an intensive research of the topic recording strong relations of the construction to properties of linear series that had not been observed at first. 

The idea of the construction of Okounkov bodies is to associate to linear series on projective varieties convex geometric objects and can be viewed as a generalization of the toric case, where to each torus-invariant divisor $D$ one an associates a polytope $P_D$. As for the toric, also in the general case Okounkov bodies are convex bodies which encode several properties of linear series as for example their volume. Although the construction of Okounkov bodies can readily be presented in abstract terms, it is notoriously hard to carry out in practice. 

A recent approach (\cite{pat}) to this problem in the case of surfaces is to find \engqq{minimal building blocks} generating all possible bodies.
In this spirit, the second author and P. \L uszcz-\'Swidecka in \cite{dp} prove that for  a smooth projective surface whose pseudo-effective cone is rational polyhedral, the Okounkov body of a big divisor with respect to a general flag decomposes as the Minkowski sum of finitely many simplices arising as Okounkov bodies of nef divisors and that finitely many such simplices suffice to construct all possible Okounkov bodies. The set of the corresponding nef divisors is called a Minkowski basis in \cite{dp}.

When generalizing the definition of a Minkowski basis to arbitrary dimension, we cannot expect Okounkov bodies of nef divisors to suffice as building blocks for all possible Okounkov bodies of big divisors, mainly because of the existence of non-nef movable divisors as we will discuss later. This leads us to the following definition. 
\begin{definition}
	Let $X$ be a projective variety of dimension $n$ and 
	$Y_{\bullet} : X = Y_{0} \supseteq Y_{1} \supseteq \dots \supseteq 
	Y_{n-1} \supseteq Y_{n} = \{pt\}$ an admissible flag on $X$.
	A collection $\set{D_1,\dots,D_r}$ of pseudo-effective divisors 
	on $X$ is called a \emph{Minkowski basis} of (Okounkov bodies on) $X$ with 
	respect to $Y_\bullet$ if
	\begin{itemize}
		\item
			for any pseudo-effective divisor $D$ on $X$ there exist non-negative numbers 
			$a_1,\dots,a_r$ such that 
			$$
				D= \sum a_iD_i \quad {\rm and} \quad 
				\Delta_{Y_\bullet}(D)=\sum a_i\Delta_{Y_\bullet}(D_i)
			$$
			\item
				the Okounkov bodies $\Delta_{Y_\bullet}(D_i)$ are indecomposable, 
				i.e., if $\Delta_{Y_\bullet}(D_i) = P_1 + P_2$ for convex bodies $P_1,P_2$, 
				then $P_j=k_j\cdot\Delta_{Y_\bullet}(D_i)$ for non-negative numbers
				$k_1,k_2$ with $k_1+k_2=1$.
	\end{itemize} 
\end{definition}

As a first step to the generalization to higher dimensional varieties of the result in \cite{dp}, in this paper we study the case of Okounkov bodies on toric varieties with respect to flags which are invariant under the action of the maximal torus. One of the main tools we are going to use comes from $\cite{lm}$, where the authors identify Okounkov bodies of a torus-invariant divisor $D$ with its moment polytope $P_D$. 

It is clear from Condition 1) that a Minkowski basis must contain a generating divisor for 
each extremal ray of the pseudo-effective cone $\Eff(X)$. 
On the other hand, we will argue that for a variety on which the section ring $\bigoplus_m H^0(X,\O_X(mD))$ of any effective divisor $D$ is finitely generated, which of course is the case for toric varieties, it suffices to consider in addition basis elements which are movable divisors, i.e., have no stable divisorial base component.

The main result of this paper can be stated as follows.
\vskip .2cm
\textbf{Theorem \ref{th:basis}.} 
{\it On a smooth projective toric variety $X$ there is a unique (up to multiplication with positive scalars) finite Minkowski basis for Okounkov bodies for all admissible torus-invariant flags. Concretely, its movable elements are given as the set
\begin{eqnarray*}
	{\rm MB_{mov}} = \left\{
						\ D\in\Div(X)\ \Bigg| \quad\begin{aligned} &\text{there is a small modification $f:X\dasharrow X'$ 		
								and}\\ 
						&\text{a nef divisor $D'$ with $D=f^\ast(D')$ such that $[D']$ is }\\
						&\text{a minimal generator of an extremal ray of $\Nef(X')$ }
					\end{aligned} \right\}.
\end{eqnarray*}
}
\vskip .2cm
Considering the toric variety $X$ as a Mori dream space, we can interpret the set ${\rm MB_{mov}}$ as the set of minimal generators of extremal rays of Mori chambers in the movable cone ${\rm Mov}(X)$ as introduced in \cite{hk}, or---in a more classical language---as minimal generators of 1-dimensional cones in the secondary fan of $X$ (\cite{gkz}, \cite{op}).
In particular, finding the Minkowski basis of $X$ is equivalent to determining the Mori chamber structure of the pseudo-effective cone, which is in turn equivalent to the description of the secondary fan. 

We demonstrate how this can be achieved by providing an effective algorithm which determines the Minkowski basis for a projective toric variety of arbitrary dimension $n$. In fact, the algorithm finds all polytopes in $\R^n$ which  arise as Okounkov bodies of movable divisors with respect to a fixed $T$-invariant flag and are indecomposable. It then follows from the proof of the theorem that this set of polytopes corresponds to the movable elements of the Minkowski basis. We give a full description---both in text and as a semi-code---of the algorithm and we prove its correctness. 

The paper is structured as follows. In Sections \ref{mink}  and \ref{toric} we recall the main definitions and properties of Okounkov bodies and toric varieties, and give a general definition of a  Minkowski basis, emphasizing the relations among the different structures. In Section \ref{alg} we prove the main theorem and discuss the algorithm.
%
Section \ref{ex} consists of examples on how the algorithm works in the case of $\P^2$ blown-up in two points and $\P^3$ blown-up in two intersecting lines. In particular, this consideration will shed light on the relevance of the existence of toric flips for the computation.

Finally, in Section \ref{decompose} we present a method to find the Minkowski decomposition of the Okounkov body for a given divisor using the output of the algorithm. This will yield the position of the class of its moving part in the movable cone.

\bigskip
{\small\noindent {\bf Acknowledgements.} } The authors would like to thank Alex K\"uronya and Tomasz Szemberg for helpful discussions and suggestions. The authors are grateful to the Pedagogical University of Cracow and IMPAN in Warsaw. We warmly
thank these institutions for their hospitality that allowed the beginning of this collaboration.

%
%
%
%
%
%
%
%
%
\section{Okounkov bodies and Minkowski decomposition}
\label{mink}

We start with a short reminder of the construction of Okounkov bodies.
Let $X$ be a normal projective variety over the complex numbers of dimension $n$. Suppose that $D$ is a Cartier divisor. To such a divisor we associate a convex compact subset of $\R^n$ called the Okounkov body. The construction requires the choice of  an admissible flag, i.e. a  sequence of subvarieties
$$Y_{\bullet} : X = Y_{0} \supseteq Y_{1} \supseteq ... \supseteq Y_{n-1} \supseteq Y_{n} = \{pt\}$$ such that 
${\rm codim } Y_{i} = i$ and $Y_{n}$ is a smooth point on each $Y_{i}$. One defines a valuation-like function
\be
	\nu : H^{0}(X, \mathcal{O}_{X}(D)) &\to& \mathbb{Z}^{n}\\
	s &\mapsto& (\nu_{1}(s), ..., \nu_{n}(s))
\ee
 in the following way: set $\nu(0)= (0,\dots,0)$ and for $0 \neq s \in H^{0}(X, \mathcal{O}_{X}(D))$ define the first coordinate of $\nu(s)$ as $\nu_{1}(s) = {\rm ord}_{Y_{1}}(s)$. After choosing a local equation of $Y_{1}$ in $X$, the section $s$ determines a section $\tilde{s_{1}} \in H^{0}(X, \mathcal{O}_{X}(D - \nu_{1}Y_{1}))$ that does not vanish identically along $Y_{1}$. Taking the restriction $s_{1} \in H^{0}(Y_{1}, \mathcal{O}_{Y_{1}}(D - \nu_{1} Y_{1}))$ we define $\nu_{2}(s) = {\rm ord}_{Y_{2}}(s_{1})$.
Continuing the procedure we obtain the vector $\nu(s) = (\nu_{1}(s), \dots, \nu_{n}(s))\in\Z^n$.
Defining the semi-group of valuation vectors as
 $$\Gamma_{Y_\bullet}(D):=  \set{ (\nu(s),m) \in \Z^n\times\N \with H^{0}(X, \mathcal{O}_{X}(mD))}
 $$ the \emph{Okounkov body} of $D$ is given by
$$\Delta_{Y_\bullet}(D) := \text{cone}_{\R^{n+1}}\left(\Gamma_{Y_\bullet}(D)\right)\cap  \left(\R^n\times \set{1}\right).$$

It is worth pointing out that the above construction can be carried out for non-complete graded linear series as well, see \cite{kk}. Okounkov bodies encode geometric and algebraic information about divisors, for instance their volume. 

Determining the Okounkov body of a given divisor is in general a quite difficult task. There are two cases in which it is possible to bypass the construction described above and to determine Okounkov bodies more directly, namely the case of surfaces and of toric varieties \cite{lm}. In the case of surfaces, Zariski decomposition of big divisors can be used to produce two functions which bound the Okounkov body. This description provides information on possible shapes of Okounkov bodies on surfaces as was worked out nicely in \cite{klm}. However, in practice it turns out that this approach is still quite involved. This fact motivated research in the direction of Minkowski decomposition. The idea stemming from \cite{pat} is to find a (preferably finite) set of divisors such that the knowledge of their Okounkov bodies suffices to determine the Okounkov bodies of arbitrary big divisors. In the case of surfaces this has been accomplished assuming the pseudo-effective cone is rational polyhedral in \cite{dp} by applying the decomposition of the big cone introduced in \cite{bks}. The result in the surface case is the following.

\begin{theorem}[\cite{dp}] Let $X$ be a smooth complex projective surface such that $\overline{{\rm Eff}(X)}$ is rational polyhedral. Fix a flag $X\supseteq L\supseteq x$, such that $L$ is a big and nef curve and $x$ is a general point of $L$. Then there exists a finite set of nef divisors ${\rm MB}(X) = \{P_{1}, \dots, P_{s}\}$ such that for every big and nef $\mathbb{Q}$-divisor $D$ there exist non-negative numbers $a_{i} \geq 0$ such that 
$$D = \sum_{i = 1}^{s} a_{i}P_{i} \qquad \text{and} \qquad \triangle(D) = \sum_{i = 1}^{s} a_{i} \triangle(P_{i}).$$
\end{theorem}
This result is constructive, i.e., its proof contains an algorithm, which allows us to determine Okounkov bodies effectively. 
\begin{example} Let $X$ be the blow up of $\mathbb{P}^2$ in two points with exceptional divisors $E_{1}, E_{2}$ and consider the flag  $Y_{\bullet}=\{l, x\}$ given by a general line and a general point on it. Then by the proof of the above result a Minkowski basis of $X$ with respect to $Y_\bullet$ is given by the set of divisors $$\{ H, H-E_{1}, H-E_{2}, 2H-E_{1}-E_{2}\},$$ where $H$ denotes the pullback of the class of a line in $\P^2$. Let us consider the divisor $D = 7H - 2E_{1} - 2E_{2}$. We can write $D= 3H + 2(2H-E_1-E_2)$ and hence obtain 
$$\triangle(D) = 3\triangle(H) + 2 \triangle(2H-E_{1} - E_{2})$$
as a Minkowski decomposition of its Okounkov body.\bigskip\\
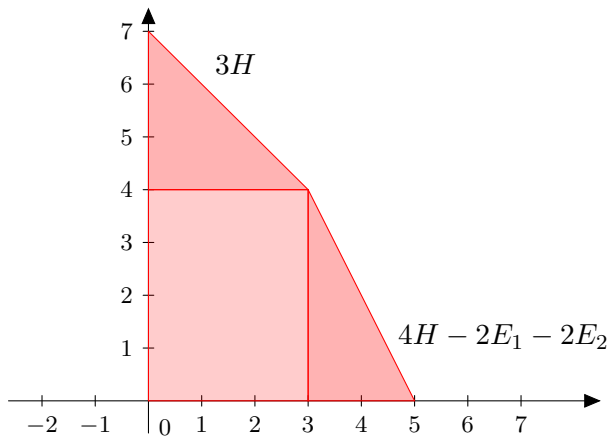
\begin{center}
\begin{tikzpicture}[line cap=round,line join=round,>=triangle 45,x=.7cm,y=.7cm]
\draw[->,color=black] (-2.63,0) -- (8.5,0);
\foreach \x in {-2,-1,1,2,3,4,5,6,7}
\draw[shift={(\x,0)},color=black] (0pt,2pt) -- (0pt,-2pt) node[below] {\footnotesize $\x$};
\draw[->,color=black] (0,-0.61) -- (0,7.45);
\foreach \y in {1,2,3,4,5,6,7}
\draw[shift={(0,\y)},color=black] (2pt,0pt) -- (-2pt,0pt) node[left] {\footnotesize $\y$};
\draw[color=black] (0pt,-10pt) node[right] {\footnotesize $0$};
\fill[color=red,fill=red,fill opacity=0.3] (0,7) -- (0,4) -- (3,4) -- cycle;
\fill[color=red,fill=red,fill opacity=0.3] (3,4) -- (3,0) -- (5,0) -- cycle;
\fill[color=red,fill=red,fill opacity=0.2] (0,4) -- (3,4) -- (3,0) -- (0,0) -- cycle;
\draw [color=red] (0,7)-- (0,4);
\draw [color=red] (0,4)-- (3,4);
\draw [color=red] (3,4)-- (0,7);
\draw [color=red] (3,4)-- (3,0);
\draw [color=red] (3,0)-- (5,0);
\draw [color=red] (5,0)-- (3,4);
\draw [color=red] (0,4)-- (3,4);
\draw [color=red] (3,4)-- (3,0);
\draw [color=red] (3,0)-- (0,0);
\draw [color=red] (0,0)-- (0,4);
\draw (1.06,6) node[anchor=south west] {$3H$};
\draw (4.5,.8) node[anchor=south west] {$4H -2E_1 -2 E_2$};
\end{tikzpicture}
\end{center}
\end{example}
In order to generalize the idea to higher dimensions, we make the following definition.
\begin{definition}
	Let $X$ be a normal projective variety of dimension $n$ and 
	$Y_{\bullet} : X = Y_{0} \supseteq Y_{1} \supseteq \dots \supseteq 
	Y_{n-1} \supseteq Y_{n} = \{pt\}$ an admissible flag on $X$.
	A collection $\set{D_1,\dots,D_r}$ of pseudo-effective $\Q$-divisors 
	on $X$ is called a \emph{Minkowski basis} of (Okounkov bodies on) $X$ with 
	respect to $Y_\bullet$ if
	\begin{itemize}
		\item
			for any pseudo-effective $\Q$-divisor $D$ on $X$ there exist 
			non-negative numbers 
			$a_1,\dots,a_r$ such that 
			$$
				D= \sum a_iD_i \quad {\rm and} \quad 
				\Delta_{Y_\bullet}(D)=\sum a_i\Delta_{Y_\bullet}(D_i)
			$$
			\item
				the Okounkov bodies $\Delta_{Y_\bullet}(D_i)$ are indecomposable, 
				i.e., if $\Delta_{Y_\bullet}(D_i) = P_1 + P_2$ for convex bodies $P_1,P_2$, 
				then $P_j=k_j\cdot\Delta_{Y_\bullet}(D_i)$ for non-negative numbers
				$k_1,k_2$ with $k_1+k_2=1$.
	\end{itemize} 
\end{definition}

As noted in the introduction, any Minkowski basis must contain divisors spanning the extremal rays of
the pseudo-effective cone $\Eff(X)$ in order to satisfy condition 1). So these elements are fixed up to multiplication with positive scalars. In case $X$ is a toric variety, or more generally a variety on which section rings $R(X,D):=\bigoplus_m H^0(X,\O_X(mD))$ of effective divisors $D$ are finitely generated, we can in fact assume the remaining basis elements to be movable. Indeed, let $D$ be any pseudo-effective divisor on $X$. We can without loss of generality assume $R(X,D)$ to be generated in degree one, so that  decomposing $D$ into its movable and fixed part $D = M+ F$, multiplication with $s_F^m$ for a defining section $s_F$ of $F$ yields isomorphisms
$$
	H^0(X,\O_X(mM))\to H^0(X,\O_X(mD))
$$
for all $m\ge1$. In particular, for any admissible flag $Y_\bullet$ we get the identity of Okounkov bodies
$$
	\Delta_{Y_\bullet}(D) = \Delta_{Y_\bullet}(M) + \Delta_{Y_\bullet}(F) = \Delta_{Y_\bullet}(M) + \nu(s_F)
$$
and $F$ lying on the boundary of $\Eff(X)$ can be written as a positive sum of basis elements spanning the extremal rays. Hence it is enough to decompose movable divisors in order to complete the  Minkowski basis.

Note that in the result on the surface case cited above 
the set MB together with classes spanning the (finitely many) 
extremal rays of $\Eff(X)$ constitutes a Minkowski basis in the above sense, 
since the Okounkov body of any pseudo-effective divisor on a surface is just 
a translate of the Okounkov body of its positive part, which 
is nef, and thus can be decomposed into bodies coming from 
the divisors in MB. Furthermore, it follows from the construction
that these divisors have either lines or triangles as Okounkov bodies, in 
particular these are indecomposable.

In the present paper we prove that for a smooth projective toric variety $X$
the movable elements of a Minkowski basis with respect to any $T$-invariant flag are 
given by the set of all movable divisors which span an extremal ray 
of the nef cone $\Nef(X')$ for some small modification of $X$. Surprisingly, this means that  
the Minkowski basis is independent of the flag. 

%
%
%
%
%

\section{Toric varieties and Okounkov bodies}
\label{toric}

In this section we collect results about toric varieties and divisors on them which we will need in the proof of the main theorem. The main references for this section are \cite{toric} and \cite{lm}.

Consider a normal projective toric variety $X=X_{\Sigma}$ corresponding to a complete fan $\Sigma$ in $N_{\R}$ (with no torus factor) with $\dim N_{\R}=n$. Recall that every $T_N$-invariant Weil divisor is represented as a sum $$D= \sum_{\rho \in \Sigma(1)} a_{\rho}D_{\rho},$$
where $\rho$ is a one-dimensional subcone (a ray), and $D_{\rho}$ is the associated $T_N$-invariant prime divisor. 
The divisor $D$ is Cartier if for every maximal dimensional subcone $\sigma \in \Sigma(n)$, the restriction $D|_{U_{\sigma}}$ is locally the divisor of a character (i.e., of the form $\text{div}(\chi^{m_{\sigma}})$ with $m_{\sigma} \in N^{\vee}=M$).

To every divisor one associates a polyhedron
$$P_D=\{m\in M_{\R} | \langle m, u_{\rho}\rangle\geq - a_{\rho} \text{ for every } \rho \in \Sigma(1) \}.$$
whose integral points represent the global sections of the divisor. We note that the toric polyhedron detects whether a given divisor is globally generated, hence nef:

\begin{theorem}[\cite{toric}{ 6.3.12}] For a divisor $D$ on a toric variety $X$ the following properties are equivalent.
\begin{enumerate}
\item $D$ is nef.
\item $\mathcal{O}_X(D)$ is generated by global sections.
\item $m_{\sigma} \in P_D$ for all $\sigma \in \Sigma(n)$.
\end{enumerate} 
\end{theorem}

\begin{example} \label{p21} Let $X$ be the blow-up of $\P^2$ in one point. The fan $\Sigma$ is given by $4$ rays, i.e., $D_1=H-E_1$, $D_2=H$, $D_3= H-E_1$ and $E_1$.  The following figure illustrates the construction of the polytope $P_D$ for the  $T$-invariant divisor $D=D_1+D_2+D_3+E_1$.
\vskip .5cm
\begin{tikzpicture}[line cap=round,line join=round,>=triangle 45,x=1.0cm,y=1.0cm]
\draw[color=black] (0,0) -- (3,0);
\draw (0.1,2.5) node [anchor=north west] {$D_2$};
\draw (2.5,0.1) node[anchor=south west] {$D_1$};
\draw (0.1, -2.5) node[anchor=west] {$E_1$};
\draw (-2.2,-1.5) node[anchor=north east] {$D_3$};
\draw[color=black] (0,-2.5) -- (0,2.5);
\draw[color=black] (0,0) -- (-2.5,-2.5);
\draw (-2,2.5) node[anchor=north east] {$\Sigma \in N$};
\end{tikzpicture}
$\qquad$
\begin{tikzpicture}[line cap=round,line join=round,>=triangle 45,x=1.0cm,y=1.0cm]
\draw[->,color=gray] (0,0) -- (3,0);
\draw[->,color=gray] (0,0) -- (0,3);
\draw[color=green!50!black] (-1.5,1.5) -- (-1.5,-1.5);
\draw[color=green!50!black] (-1.5,-1.5) -- (3,-1.5);
\draw[color=green!50!black] (3,-1.5) -- (0,1.5);
\draw[color=green!50!black] (0,1.5) -- (-1.5,1.5);
\draw (0, 1.5) node [color=red]{\LARGE{\textbullet}};
\fill[color=red,fill=red,fill opacity=0.3] (-1.4,1.4) -- (-1.4,-1.4) -- (2.8,-1.4) -- (0, 1.4 ) -- cycle;
\draw (1, 1) node[anchor=west] {$P_D$};
\draw (1, -2.5) node[anchor=west] {\emph{}};
\draw (-1,3) node[anchor=north east] {$M=N^{\vee}$};
\end{tikzpicture}
\end{example}

\begin{theorem} [\cite{toric}{ 6.2.8}] \label{main} Let $D$ be a basepoint 
free Cartier divisor on a complete toric variety, and let 
$X_D$ be the toric variety of the polytope $P_D \subseteq M_{\R}$. 
Then any refinement  $\Sigma$ of the fan $\Sigma_{P_D}$ induces a 
proper toric morphism 
$$
	\varphi : X_{\Sigma} \to X_{P_D}.
$$
Furthermore, $D$ is linearly equivalent to the pullback via 
$\varphi$ of the ample divisor on $X_{P_D}$ coming from $P_D$.
\end{theorem}

\begin{remark}
Theorem \ref{main} implies  in particular  that the polytope of the pullback 
of a $T$-invariant divisor  by a toric morphism  is isomorphic to the original polytope. \bigskip\\
\end{remark}
The following statements also are important consequences of \ref{main}.

\begin{proposition} [\cite{toric}{ 6.2.13}]  \label{ref} 
	Let $P$ and $Q$ be lattice polytopes in $M_{\R}$. Then:
	\begin{enumerate}
		\item $Q$ is an $\N$-Minkowski summand of $P$ if and only if 
			$\Sigma_P$ refines $\Sigma_Q$.
		\item 
			$\Sigma_{P+Q}$ is the coarsest common refinement of $\Sigma_P$ 
			and $\Sigma_Q$.
	\end{enumerate} 
\end{proposition}

\begin{corollary} [\cite{toric}{ 6.2.15}]
	Let $P$ be a full dimensional 
	lattice polytope in $M_{\R}$. Then a
	polytope $Q \subseteq M_{R}$ is an $\N$-Minkowski summand of $P$ if 
	and only if there is a torus invariant
	basepoint free Cartier divisor $A$ on $X_P$ such that $Q = P_A$.
\end{corollary}

\begin{lemma} [\cite{fulton}{ Section 3.4}] \label{add} If $D$ and $E$ are nef $T$-invariant divisors, then 
$$P_{D+E}= P_E+P_D.$$
\end{lemma}
Finally, the following statement links the polytope $P_D$ to the Okounkov body $\Delta(D)$.

\begin{proposition}[\cite{lm}{ 6.1}]  \label{lm} Let $X_\Sigma$ be a smooth projective toric variety, and let $Y_{\bullet}$ be a flag given as a complete intersection of a set of $T$-invariant divisors generating a maximal cone $\sigma\in\Sigma$. Given any big line bundle $\mathcal{O}_X(D)$ on X, such that $D|_{\mathcal{U}_{\sigma}} = 0$, then
$$ \Delta(\mathcal{O}_X(D)) = \varphi(P_D),$$
where $\varphi$ is an $\R$-linear map.
\end{proposition}


\begin{example}  In the situation of Example \ref{p21} consider the $T$-invariant flag $Y_{\bullet}=\{D_3, D_3\cap E_1\}$. The element $ D'= 2D_2 +D_1 \in |D|$ satisfies the condition $D'|_{\mathcal{U}_{\sigma}} = 0$. Thus, the Okounkov body $\Delta(D)$ with respect to $Y_{\bullet}$ is given by the image of $P_{D'}$ under the linear transformation $\varphi$ which maps the rays spanned by $(-1,0)$ and $(1,-1)$ to the coordinate axes.
\vskip .5cm
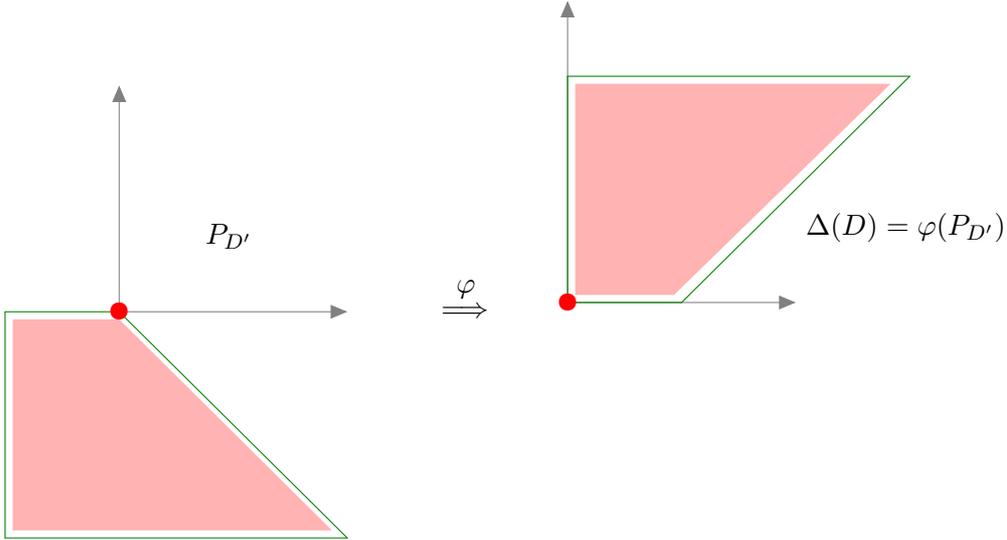
\begin{tikzpicture}[line cap=round,line join=round,>=triangle 45,x=1.0cm,y=1.0cm]
\draw[->,color=gray] (0,0) -- (3,0);
\draw[->,color=gray] (0,0) -- (0,3);
\draw[color=green!50!black] (-1.5,0) -- (-1.5,-3);
\draw[color=green!50!black] (-1.5,-3) -- (3,-3);
\draw[color=green!50!black] (3,-3) -- (0,0);
\draw[color=green!50!black] (0,0) -- (-1.5,0);
\draw (0, 0) node [color=red]{\LARGE{\textbullet}};
\fill[color=red,fill=red,fill opacity=0.3] (-1.4,-0.1) -- (-1.4,-2.9) -- (2.8,-2.9) -- (0, -0.1 ) -- cycle;
\draw (1, 1) node[anchor=west] {$P_{D'}$};
\draw (4.1,0) node[anchor=west] {$\Longrightarrow$};
\draw (4.3,.3) node[anchor=west] {$\varphi$};
\draw (-1,3) node[anchor=north east] {\emph{}};
\end{tikzpicture}
$\quad$
\begin{tikzpicture}[line cap=round,line join=round,>=triangle 45,x=1.0cm,y=1.0cm]
\draw[->,color=gray] (0,0) -- (3,0);
\draw[->,color=gray] (0,0) -- (0,4);
\draw[color=green!50!black] (1.5,0) -- (4.5,3);
\draw[color=green!50!black] (4.5,3) -- (0,3);
\draw[color=green!50!black] (0,3) -- (0,0);
\draw[color=green!50!black] (0,0) -- (1.5,0);
\draw (0, 0) node [color=red]{\LARGE{\textbullet}};
\fill[color=red,fill=red,fill opacity=0.3] (1.4,0.1) -- (4.25,2.9) -- (0.1,2.9) -- (0.1, 0.1 ) -- cycle;
\draw (3, 1) node[anchor=west] {$\Delta(D)= \varphi(P_{D'})$};
\draw (1, -3) node[anchor=west] {\emph{}};
\end{tikzpicture}
\end{example}

%
%
%

\section{Minkowski bases for toric varieties}
\label{alg}
In this section we characterize movable Minkowski basis elements 
for a smooth projective toric variety of arbitrary dimension and 
present an algorithm to obtain them, given as data the complete 
fan defining the toric variety.

\subsection{Characterization of Minkowski basis elements}
Let $X$ be a smooth projective toric variety of dimension $n$. 
\begin{remark}
Note that a Minkowski basis for any given $T$-invariant flag 
will also be a Minkowski basis for any other $T$-invariant flag. 
This is due to Theorem \ref{lm} which says that 
the (transpose of the) Okounkov body of a big divisor $D$ on $X$ with respect to 
a flag $Y_\bullet$ is the linear transform of the polytope $P_D$.
Now, Minkowski summation obviously respects linear transformation, 
so both properties of a Minkowski basis can be checked immediately 
on the level of $T$-invariant polytopes, which is what we will always do 
in the rest of the paper.\bigskip\\
\end{remark}
Our aim is to prove the following.
\begin{theorem}\label{th:basis}
		On a smooth projective toric variety $X$ equipped with an admissible 
		$T$-inariant flag $Y_\bullet$  there is a unique (up to multiplication with positive 
		scalars) finite Minkowski basis for Okounkov bodies. Concretely, its movable elements are given as the set
	\begin{eqnarray*}
		{\rm MB_{mov}} = \left\{
							\ D\in\Div(X)\ \Bigg| \quad\begin{aligned} &\text{there is a small 
							modification $f:X\dasharrow X'$ 		
									and}\\ 
							&\text{a nef divisor $D'$ with $D=f^\ast(D')$ such that $[D']$ is }\\
							&\text{a minimal generator of an extremal ray of $\Nef(X')$ }
						\end{aligned} \right\}.
	\end{eqnarray*}
\end{theorem}

\begin{remark}The theorem can be rephrased in the following way: finding the
movable elements of a Minkowski
basis for a toric variety is equivalent to finding the rays in the 
secondary fan contained in the moving cone. (For the definition and 
an analysis of the secondary fan see for example \cite[14.4]{toric}).\bigskip\\
\end{remark}
Let us first prove that the set of divisors spanning extremal rays of the nef cone on some
small modification satisfies the second condition for a Minkowski basis. 
\begin{proposition} \label{second}
	For a $T$-invariant movable divisor $D$ the polytope $P_D$ is indecomposable 
	if and only if there exists a small modification $f: X\dasharrow X'$ and a divisor $D'$ 
	spanning an extremal ray of $\Nef(X')$ with $D=f^\ast(D')$.
\end{proposition}
\begin{remark}
	This result is slightly stronger than what we need in order to prove the
	theorem. In fact, it shows that a Minkowski basis 
	is unique up to multiplication and its movable elements correspond 
	exactly to the set of rays in the secondary fan. Therefore, finding 
	all indecomposable polytopes coming from movable divisors recovers this Minkowski basis. This is 
	exactly what the algorithm described in the next section does.	
\end{remark}

\begin{proof} 
First note that we can assume $D$ to be nef on $X$ since being movable, it is 
the pullback of a nef divisor under a small modification $f:X\dasharrow X'$ and 
the polytopes agree since $f$ does not alter the rays in the fan defining 
$X$, but only changes higher dimensional cones. 

Let us consider $X_{P_D}$, the variety given by the normal fan of $P_D$. 
Note that $\Sigma$ is a refinement of $\Sigma_{P_D}$, in particular we 
obtain a proper toric morphism $\varphi: X \to X_{P_D}$.

Let us now suppose that $D$ spans an extremal  ray of ${\rm Nef}(X)$ 
and that $P_D=Q+N$ is a Minkowski decomposition. Since $P_D$ is a full dimensional polytope with 
respect to $X_{P_D}$,  there exist nontrivial basepoint 
free T-invariant, hence nef, Cartier divisors $A$ and $B$ on $X_{P_D}$ 
such that $Q=P_A$ and $N=P_B$.

But now we have that $\varphi^*(A)$, $\varphi^*(B)$ are nef divisors 
on $X$ and since toric polytopes remain invariant under pullback, we have
$$P_D= P_{\varphi^*(A)}+ P_{\varphi^*(B)} \quad \mbox{and} \quad D = 
\varphi^*(A)+ \varphi^*(B)$$
hence $\varphi^*(A)$ and $\varphi^*(B)$ are equivalent to multiples of $D$ since $[D]$ 
spans an extremal ray and thus the decomposition $P=Q+N$ is trivial.

For the opposite implication, assume $D$ does not lie in an extremal ray
of $\Nef(X)$. We can then find nef divisors $A$ and $B$ which do not lie in the 
ray spanned by $[D]$ such that $D=A+B$. Now for globally 
generated divisors we know that the polytope 
of a sum coincides with the Minkowski sum of the polytopes, i.e.,
$$
	P_D=P_A+P_B,
$$
showing that $P_D$ is decomposable in a non-trivial way.
\end{proof}

\begin{proof}[Proof of the theorem]
	It remains to show the first property.  
	Let therefore $\sigma\in\Sigma(n)$ be the  maximal cone defining the flag $Y_\bullet$.
	We consider $X$ as a Mori dream space. Let $D$ be a $T$-invariant pseudo-effective 
	divisor on $X$. Its class lies in some Mori chamber $\Sigma_g$ for a contracting rational
	map $g: X\dasharrow Y$. 
	In particular by \cite[Proposition 1.11]{hk} $D$ has a decomposition $D=P+N$
	 into its movable and fixed parts such that $N$ is supported on the set of exceptional 
	divisors of $g$ and $D$ is the pullback of a nef divisor on $Y$. 
	
	Fixing a defining section 
	$s_N$ for $N$, each section $s\in H^0(X,\O_X(mD))$ factors as $s=s'\cdot s_N^m$ for a 
	section $s'\in H^0(X,\O_X(mP))$. Therefore, the Okounkov
	body of $D$ is a translate by $\nu(s_N)$ of the body of $P$, which is a 
	movable divisor. Without loss of generality we can pick $P$ as a representative
	of its numerical class which	satisfies the condition $P\mid_{U_\sigma}=0$, and
	hence can identify $\Delta_{Y_\bullet}$ with the moment polytope $P_{P}$ via the 
	linear transform $\varphi_\sigma$. 
	By  \cite[Proposition 1.11]{hk} there is a small toric modification $f:X\dasharrow X'$ 
	such that $P$ is the pullback of a nef divisor $P'$ on $X'$.
	 Note that such a modification does not
	alter the polytopes, i.e., $P_P=P_{P'}$ since their construction only depends on the 
	rays in the respective fans, which are unchanged by the small modification.

	The polytope corresponding to $P'$ can be decomposed as a Minkowski sum 
	of polytopes of divisors, say $D_1',\ldots,D_r'$, spanning extremal rays of $\Nef(X')$. 
	This immediately follows from the fact that we can write $P'$ as a 
	non-negative $\Q$-linear combination 
	$$
		P'=\sum_{i=1}^r a_iD_i'
	$$
	together with Lemma \ref{add}, which states the 
	additivity of moment polytopes for nef divisors, i.e.,
	$$
		P_{P'}=\sum_{i=1}^r a_i P_{D_i'}.
	$$
	The nef divisors $D_1',\ldots,D_r'$ correspond to movable divisors $D_1,\ldots,D_r$ 
	on $X$ and the corresponding polytopes agree, again since the modification does not
	change the rays in the fan defining $X$ and $X'$, respectively.
	
	By passing to representatives of the $D_i$ which satisfy $D_i\mid_{U_\sigma}=0$, 
	$\varphi^{-1}$ maps the polytopes $P_{D_i}$ to the Okounkov bodies 
	$\Delta_{Y_\bullet}(D_i)$ and respects Minkowski sums, i.e.,
	$$
		\Delta_{Y_\bullet}(P) =\sum a_i\Delta_{Y_\bullet}(D_i).
	$$

	The proof of the Theorem then follows from Proposition \ref{second}.
\end{proof}

\begin{remark}
	Note that in the above argument we do not consider the Okounkov body 
	of $D'$ on the flipped model $X'$ but only the polytope $P_{D'}$, which we compare to 
	$P_D$. Therefore, we do not apply Proposition \ref{lm} to $X'$, and in particular we 
	need not assume $X'$ to be smooth. By the same token, we do not need additional 
	assumptions on $T$-invariant flags on $X'$.
\end{remark}

\subsection{How to find a Minkowski basis}
Let us consider a toric variety $X=X_\Sigma$, $\mbox{dim}(X)=n$ and let $d= \#(\Sigma(1))-n$, where as usual $\Sigma(1)$ denotes the set of rays in  $\Sigma$. We will first informally describe the idea of the algorithm determining the movable elements of the Minkowski basis  of $X_{\Sigma}$. A formal description of the algorithm follows afterwards. In the following, $R_m$ will denote a set of rays of the fan $\Sigma$ and $R^*_m$ denotes the set of half-spaces dual to given rays in $R_m$. Furthermore, $P_m$ will be the set of points corresponding to the vertices of $R^*_m$. By $\Delta_i=\text{convex hull}(P_m)$ denote the convex hull generated by the points in $P_m$.

We also define the half-spaces $H_t(\rho_i)=\{x\in N_{\R} | \langle x, \rho_i \rangle \geq -t\}$ and 
the bounding hyperplane $\bar{H}_t(\rho_i) =  \{x\in N_{\R} | \langle x, \rho_i \rangle = -t\}$.

The idea behind the algorithm is quite straightforward: every \engqq{slope} appearing in some toric polyhedron has to occur in the polytope of at least one of the Minkowski basis elements. The steps of the algorithm are the following:
\begin{itemize}
\item Fix a $T$-invariant flag, i.e., choose a set of $n=\dim(X)$ rays that generate one of the cones $\sigma \in \Sigma(n)$. This cone corresponds to a cone $\sigma^{\vee}$ in the dual space $M_{\R}$. Denote by $\Sigma_0$ the (non-complete) fan  generated by the rays of $\sigma$.

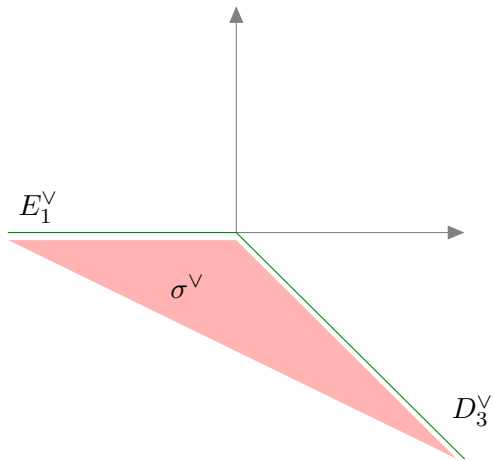
\qquad\qquad\quad\begin{tikzpicture}[line cap=round,line join=round,>=triangle 45,x=1.0cm,y=1.0cm]
\draw[->,color=gray] (0,0) -- (3,0);
\draw[->,color=gray] (0,0) -- (0,3);
\draw[color=green!50!black] (0,0) -- (-3,0);
\draw (-3,0) node [anchor=south west] {$E_1^{\vee}$};
\draw[color=green!50!black] (0,0) -- (3,-3);
\draw (2.7,-2.7) node[anchor=south west] {$D_3^{\vee}$};
\fill[color=red,fill=red,fill opacity=0.3] (0,-0.1) -- (-3,-0.1) -- (2.9,-3) -- cycle;
\draw (-1,-1) node[anchor=south west] {$\sigma^{\vee}$};
\draw (1, -3.2) node[anchor=north west] {\emph{}};
\end{tikzpicture}
\vskip .5cm
\item In the first step, fix a slope, i.e., a ray $\rho_1 \in \Sigma(1)\backslash \sigma(1)$. To this ray consider the corresponding an hyperplane $H_1(\rho_1)$. 

\item $H_1(\rho_1)$ intersects an edge of $\sigma^{\vee}$ if and only if the corresponding hyperplanes generate a convex cone in $N_{\R}$. If it intersects all the rays, then the intersection of the half-space and the cone $\sigma^\vee$ is a simplex, hence corresponds to a Minkowski basis element. If it does not intersect all rays, it is unbounded. Either way, call this polyhedron $\Delta_1$ and the corresponding fan $\Sigma_1$.

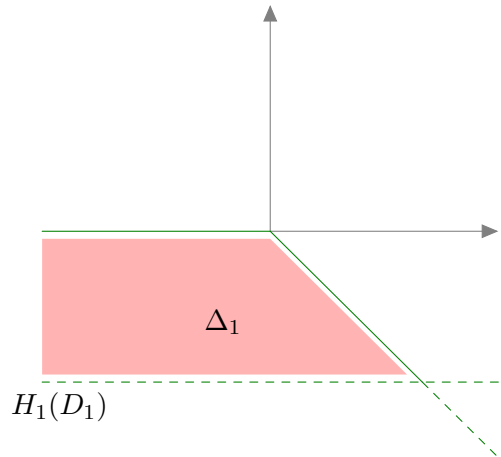
\qquad\qquad\begin{tikzpicture}[line cap=round,line join=round,>=triangle 45,x=1.0cm,y=1.0cm]
\draw[->,color=gray] (0,0) -- (3,0);
\draw[->,color=gray] (0,0) -- (0,3);
\draw[dashed, color= green!50!black] (3,-2) -- (-3,-2);
\draw[color=green!50!black] (0,0) -- (-3,0);
\draw[dashed,color=green!50!black] (2,-2) -- (3,-3);
\draw (-2,-2) node [anchor=north east] {$H_1(D_1)$};
\draw[color=green!50!black] (0,0) -- (2,-2);
\fill[color=red,fill=red,fill opacity=0.3] (0,-0.1) -- (1.8,-1.9) -- (-3,-1.9) -- (-3, -0.1) -- cycle;
\draw (-1,-1.5) node [anchor=south west] {$\Delta_1$};
\end{tikzpicture}

\item If $\Delta_1$ is not a polytope, add an additional ray, say $\rho_2$, into $\Sigma_1$, i.e., intersect the polyhedron $\Delta_1$ with an additional half-space. There are two different ways to intersect while preserving indecomposability: either take $\Delta_2$ to be the intersection of $\Delta_1$ with the half-space $H_t(\rho_2)$, where $t$ is the minimal positive number such that all vertices of $\Delta_1$ are contained in $H_t(\rho_2)$, or with the half-space $H_0(\rho_2)$. Apply both intersections consecutively.

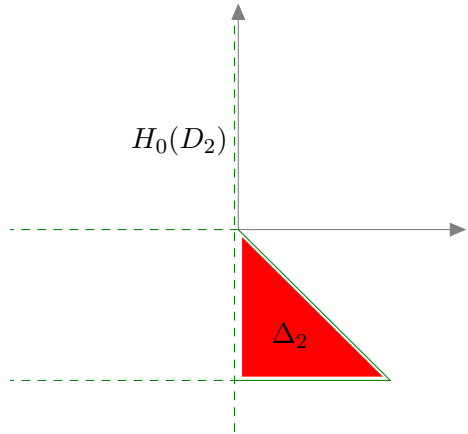
\qquad\qquad\quad\begin{tikzpicture}[line cap=round,line join=round,>=triangle 45,x=1.0cm,y=1.0cm]
\draw[->,color=gray] (0,0) -- (3,0);
\draw[->,color=gray] (0,0) -- (0,3);
\draw[color= green!50!black] (2,-2) -- (0,-2);
\draw[dashed, color= green!50!black] (0,-2) -- (-3, -2);
\draw[dashed, color=green!50!black] (0,0) -- (-3,0);
\draw[dashed,color=green!50!black] (-0.05,2.7) -- (-0.05,-2.7);
\draw (0,1.5) node [anchor=north east] {$H_0(D_2)$};

\draw[color=green!50!black] (0,0) -- (2,-2);
\fill[color=red,red] (0.05,-0.1) -- (1.9,-1.95) -- (0.05,-1.95) -- cycle;
\draw (0.3,-1.7) node [anchor=south west] {$\Delta_2$};
\draw (0,-3) node [anchor=south east] {\emph{}};
\end{tikzpicture}

\item
Keep intersecting until all $\rho\in\Sigma_1\backslash \sigma(1)$ are exhausted. Denote by $\Delta_k$ the region bounded by the already existing hyperplanes after $k$ steps. By Proposition \ref{prop:indecomp} $\Delta_k$ is indecomposable if it is bounded. Check whether it is the polytope of a divisor on $X$. If this is the case this divisor is a Minkowski base element. Start again with a different $\rho_1$.
\end{itemize}

%
%
%
%
			

We now give a formal description of the above sketch. 

\pagebreak
\begin{algorithm}\label{algo}\rm  Algorithm \emph{TMB}\\
{
{\footnotesize
   \begin{lines}

     \fbox{\parbox{\linewidth}{
      Input: fan $\Sigma$ with $n+d$ rays
      Output: Minkowski basis for $X_\Sigma$
      Variables:
		      \+$R$, an array of length $d$; each entry is a set of rays: 
		      \+\+ records the rays yet to use.
      		\+$M$, an array of length $d$; each entry is either a ray or empty:
      		\+\+ records rays currently use. 
      		\+$P$, an array of length $d$; each entry is a set of points:
      		\+\+ records points left to put hyperplanes through.
      		\+$\Delta$, an array of length $d$; each entry is a polyhedron:
      		\+\+ records the resulting polyhedron in the current step.
%
\begin{multicols}{2}
   for k from 1 to $d$ do
   \+$R_k \gets \Sigma(1)$;
   \+$M_k,\Delta_k \gets \emptyset$;
   \+$P_k\gets \emptyset$;
   end do;
   for $l$ from $1$ to $d$ do
   \+ $R_1\gets \Sigma(1)\backslash \rho_l;$
   \+ $\Delta_1\gets H_1(\rho_l)\cap \sigma^\vee;$
   \+ if (\emph{CorrespondsToDivisor}($\Delta_{1}$)) then
      \+\+ $TMB\gets TMB \cup \set{\Delta_{1}}$;
      \+end if;
      \+ $P_2\gets \emph{Vertices}(\Delta_1)$;
   \+while ($R_2\neq\emptyset$) do
      \+\+ $m\gets \max\set{n\with P_n\neq \emptyset}$;
      \+\+ if ($M_m = \emptyset$) then
      	\+\+\+ Pick $\rho \in R_m$;
      	\+\+\+$R_m\gets R_{m}\backslash \rho$;
      	\+\+\+ for $k$ from $m+1$ to $d$ do  
      		\+\+\+\+ $R_k\gets R_{m-1}\backslash \rho$;
      		\+\+\+\+ $M_k\gets \emptyset$;
      	\+\+\+ end do;
      	\+\+\+ $P_m\gets $ \emph{Vertices}$(\Delta_{m-1})$;
      	\+\+\+ $M_m\gets \rho$;
      \+\+end if;
      
      \+\+ $t\gets \min_{p\in P_m}\left(\emph{Solve$_s$}(p\in \overline H_s(\rho))\right)$;
      
      \+\+ if $(t>0)$ then
      
      \+\+\+ $\Delta_{m} \gets \Delta_{m-1}\cap H_t(\rho)$;
      \+\+\+ if (\emph{CorrespondsToDivisor}($\Delta_{m}$)) then
      \+\+\+\+ $TMB\gets TMB \cup \set{\Delta_{m}}$;
      \+\+\+ end if;
      \+\+\+ $P_m\gets \set{0}$;

      \+\+ else 
      
      \+\+\+ $\Delta_{m} \gets \Delta_{m-1}\cap H_0(\rho)$;
      \+\+\+ if (\emph{CorrespondsToDivisor}($\Delta_{m}$)) then
      \+\+\+\+ $TMB\gets TMB \cup \set{\Delta_{m}}$;
      \+\+\+ end if;
      \+\+\+ $s\gets \max\set{n\le m-1\with R_n\neq \emptyset}$;
      \+\+\+ if $(m<s)$ then
      \+\+\+\+ $P_m\gets \emptyset$;
      \+\+\+\+ $P_{m+1} \gets  \emph{Vertices}(\Delta_{m})$;   
      \+\+\+\+ $M_m\gets \emptyset$;
      \+\+\+ else
      \+\+\+\+ $M_m\gets \emptyset$;
      \+\+\+\+ $P_m\gets \emptyset$;
      \+\+\+\+ if ($P_s=\emptyset)$ then
      \+\+\+\+\+ $P_s\gets \emph{Vertices}(\Delta_{s-1})$;     
      \+\+\+\+ end if;
      \+\+\+ end if;
      
      \+\+ end if;
      \+end do;
    end do;  

return TMB;
\end{multicols}
}}
   \end{lines}
}}

\end{algorithm}

\subsection{How the algorithm works}

In this subsection we will prove correctness of the algorithm. Note that it clearly terminates since only finitely many combinations of rays need to be checked. 

\begin{proposition}\label{prop:indecomp} Every polytope $\Delta$ in the output of the algorithm is indecomposable.
\end{proposition}

\begin{proof} 
We prove inductively that in each step of the algorithm the fan $\Sigma_k$ remains minimal in the sense that 
it is not the refinement of any other fan with the same convex span. In particular this means that as soon $\Delta_k$ is bounded, it is indecomposable by Proposition \ref{ref}.\\
$\Sigma_0$ is clearly minimal.
Let us suppose that $\Sigma_k$ for $k\ge0$ is a minimal fan.
When adding a face to the dual polyhedron $\Delta_k$ that is tangent to a vertex of $\Delta_k$, this is equivalent to:
\begin{itemize}
\item cancel in the fan all the cones corresponding to faces completely contained in the complement of the half-space, so that the resulting fan is minimal again
\item cancel in the fan the cone corresponding to the vertex (and all of its faces but the rays)
\item add a ray that is not contained in the convex span of the rays we are left with
\item construct all the possible cones with the new ray and the rays of the last cone we canceled
\item  if there is a containment of cones, then the ray of the inner one that is properly contained gets canceled from the new fan (we do not admit star subdivisions)
\end{itemize}

Is the new fan minimal? Yes, in fact: 
\begin{itemize}
\item by construction we cannot cancel any ray for the new cones we have, since we already canceled all the star subdivisions

 \item none of the other rays can get canceled by minimality of the fan at the previous stage.
\end{itemize}
%
\end{proof}

\begin{remark} Note that not every  indecomposable polytope of $T$-invariant divisors is 
found by the algorithm for a fixed maximal cone $\sigma\in\Sigma(n)$, but only those 
with a vertex in the origin. In particular, the output contains those irreducible polytopes 
which come from divisors whose base locus does not contain the $T$-invariant subvarieties coming from the cone $\sigma$. 
In case the variety admits a flip it is possible that not all the indecomposable polytopes corresponding to extremal rays of the movable cone will have the origin as a vertex for a fixed flag. To overcome this issue, it will be enough to repeat the algorithm for flags on toric modifications of $X$, given by a  maximal cones, and to compare the extra elements arising with the basis found in the previous steps.
\end{remark}
 
\begin{proposition} Let $X=X_{\Sigma}$ be a toric variety and let $Y_{\bullet}$ be a flag given by the complete intersection of the generators of a maximal subcone. Let $D$ be a $T$-invariant divisor such that no element of the chosen flag is contained in its base locus. Then if $P_D$  is indecomposable then $D$ is in the output of the algorithm.
\end{proposition}
\begin{proof} 
This is given by the construction of the polytopes in the algorithm. 

Denote by $\sigma$ the cone in $M_\R$ corresponding to the flag $Y_\bullet$.
Without loss of generality, we can assume $D$ to be a member of its linear series with
$D_{|\sigma}=0$. This guarantees $P_D$ to be contained in the dual cone $\sigma^\vee$. 
Write $D=\sum a_iD_i$ where the $D_i$'s are the prime $T$-invariant divisors corresponding to the rays $\rho_1,\dots,\rho_r$ of $\Sigma_X$. We can assume by reordering that the first $s$ of the $D_i$'s are those which correspond to facets of $P_D$. We describe how the algorithm finds $P_D$.

Let us define $\Delta_1$  to be $\sigma^{\vee}$. As the first step we intersect $\Delta_1$ with all the half-spaces $H_0(\rho_i)$ where $P_D$ is contained in $\overline H_0(\rho_i)$, say for all $s+1\le i\le s'$ to obtain $\Delta_2$. Note that $P_D$ is an indecomposable full dimensional lattice polytope in the intersection of $N_\R$ with these half-spaces, thus in the dual lattice $M'$ of this intersection its normal fan is not the refinement of any complete fan by Proposition \ref{ref}.. 

Let now $\Delta_3$ be the intersection of $\Delta_2$ with $H_{a_j}(\rho_j)$ for which $a_j=0$ and $j \leq s$. This will give the smallest cone containing $P_D$.


Since the fan is convex and complete, we can choose a ray $\rho_1$ that is contained in a cone adjacent to $\Delta_3^{\vee}$ and such that the induced completion in the span of the fan is contained in $\Sigma_{P_D}$. Up to rescaling we can assume that $a_1=1$. This will be the starting point for the iterative part of the algorithm. Let us define $P_1=\{\mbox{vertices of } \bar H_1(D_1)\cap \Delta_3\}$. Then, choosing the ray $\rho_2$ with the same criterion, if $t= \inf\{s| H_s(\rho_2) \supseteq P_1\}$, then by construction $t=a_2$. Iterating this step with increasing the number of intersection points each time the polytope will be fully reconstructed.

\end{proof}
%

%
%
%
%
%
%
%

\section{Examples}
\label{ex}
\subsection{The blow-up of $\mathbb{P}^2$ in two points}

Let us consider the blow-up op $\mathbb{P}^2$ in two points. The fan $\Sigma$ is given by $5$ rays, corresponding to divisor classes $D_1= H-E_1, \; D_2= H-E_2, \; D_3= H-E_1 - E_2, \; E_1 \mbox{ and } E_2$. Let us fix the flag given by $Y_{\bullet}=\{D_3, D_3\cap E_1\}$.
\vskip .5cm
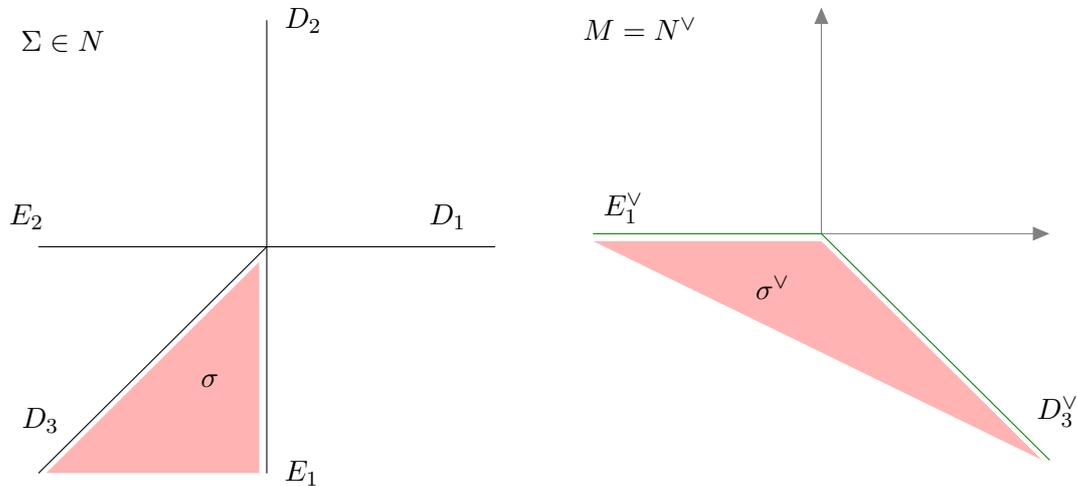
\begin{tikzpicture}[line cap=round,line join=round,>=triangle 45,x=1.0cm,y=1.0cm]
\draw[color=black] (-3,0) -- (3,0);
\draw (0.1,3) node [anchor=west] {$D_2$};
\draw (2,0.1) node[anchor=south west] {$D_1$};
\draw (-2.8,0.1) node[anchor=south east] {$E_2$};
\draw (0.1, -3) node[anchor=west] {$E_1$};
\draw (-2.6,-2.6) node[anchor=south east] {$D_3$};
\draw[color=black] (0,-3) -- (0,3);
\draw[color=black] (0,0) -- (-3,-3);
\fill[color=red,fill=red,fill opacity=0.3] (-0.1,-0.2) -- (-2.9,-3) -- (-0.1,-3) -- cycle;
\draw (-2,3) node[anchor=north east] {$\Sigma \in N$};
\draw (-1,-2) node[anchor=south west] {$\sigma$};
\end{tikzpicture}
$\qquad$
\begin{tikzpicture}[line cap=round,line join=round,>=triangle 45,x=1.0cm,y=1.0cm]
\draw[->,color=gray] (0,0) -- (3,0);
\draw[->,color=gray] (0,0) -- (0,3);
\draw[color=green!50!black] (0,0) -- (-3,0);
\draw (-3,0) node [anchor=south west] {$E_1^{\vee}$};
\draw[color=green!50!black] (0,0) -- (3,-3);
\draw (2.7,-2.7) node[anchor=south west] {$D_3^{\vee}$};
\fill[color=red,fill=red,fill opacity=0.3] (0,-0.1) -- (-3,-0.1) -- (2.9,-3) -- cycle;
\draw (-1,-1) node[anchor=south west] {$\sigma^{\vee}$};
\draw (-1.5,3) node[anchor=north east] {$M=N^{\vee}$};
\draw (1, -3.2) node[anchor=north west] {\emph{}};
\end{tikzpicture}
\vskip .5cm

We can now run the algorithm. In the following table it is recorded what is saved in each step in $R$, $M$, $P$ and $D$. 

\renewcommand{\arraystretch}{1.3}
\begin{table}[ht]\centering\footnotesize
      \begin{tabular}[t]{c|c|c|c}\tline
         $R $ &  $M$   &$P$ &       D \\\tline
 	$[[D_1,D_2, E_2],[D_1, D_2,E_2],[D_1, D_2,E_2]]$   &   $[[\emptyset],[\emptyset],[\emptyset]]$   & $[\emptyset]$ &        \\ 	

	
	$[[D_2, E_2],[D_2,E_2],[D_2,E_2]]$   &   $[[D_1],[\emptyset],[\emptyset]]$   & $[[\emptyset],[(0,0), (-1,0)],[\emptyset]]$ &              $ -$\\ 	
	
	$[[D_2, E_2],[E_2],[E_2]]$   &   $[[D_1],[D_2],[\emptyset]]$   & $[[\emptyset],[(0,0)],[(0,0),(-1,0)]]$ &         $ -$  (1)\\

	$[[D_2, E_2],[E_2],[\emptyset]]$   &   $[[D_1],[D_2],[E_2]]$   & $[[\emptyset],[(0,0)],[(0,0)]]$ &        $D_1$  \\
	$[[D_2, E_2],[E_2],[\emptyset]]$   &   $[[D_1],[D_2],[E_2]]$   & $[[\emptyset],[(0,0)],[\emptyset]]$ &        $D_1$  \\

	$[[D_2, E_2],[E_2],[E_2]]$   &   $[[D_1],[D_2],[\emptyset]]$   & $[[\emptyset],[\emptyset],[(0,0),(-1,0)]]$ &         $ -$  \\

	$[[D_2, E_2],[E_2],[\emptyset]]$   &   $[[D_1],[D_2],[E_2]]$   & $[[\emptyset],[\emptyset],[(0,0)]]$ &        $D_1$  \\
	$[[D_2, E_2],[E_2],[\emptyset]]$   &   $[[D_1],[D_2],[E_2]]$   & $[[\emptyset],[\emptyset],[\emptyset]]$ &        $D_1$  \\

	$[[D_2, E_2],[D_2],[D_2]]$   &   $[[D_1],[E_2],[\emptyset]]$   & $[[\emptyset],[(0,0)],[(0,0)]]$ &     $-$  (2)\\

	$[[D_2, E_2],[D_2],[\emptyset]]$   &   $[[D_1],[E_2],[D_2]]$   & $[[\emptyset],[(0,0)],[(0,0)]]$ &       $D_1$ (3)  \\
	$[[D_2, E_2],[D_2],[\emptyset]]$   &   $[[D_1],[E_2],[D_2]]$   & $[[\emptyset],[(0,0)],[\emptyset]]$ &       $D_1$  \\

	$[[D_2, E_2],[D_2],[D_2]]$   &   $[[D_1],[E_2],[\emptyset]]$   & $[[\emptyset],[\emptyset],[(0,0)]]$ &     $-$  \\

	$[[D_2, E_2],[D_2],[\emptyset]]$   &   $[[D_1],[E_2],[D_2]]$   & $[[\emptyset],[\emptyset],[(0,0)]]$ &       $D_1$  \\
	$[[D_2, E_2],[D_2],[\emptyset]]$   &   $[[D_1],[E_2],[D_2]]$   & $[[\emptyset],[\emptyset],[\emptyset]]$ &       $D_1$  \\


$[[D_1, E_2],[D_1,E_2],[D_1,E_2]]$   &   $[[D_2],[\emptyset],[\emptyset]]$   & $[[\emptyset],[(0,0), (1,-1)],[\emptyset]]$ &                $-$ (4)\\ 	


$[[D_1, E_2],[E_2],[E_2]]$   &   $[[D_2],[D_1],[\emptyset]]$   & $[[\emptyset], [(0,0)],[(0,0), (1, -1)]$ &              $-$ (5)\\ 	


$[[D_1, E_2],[E_2],[\emptyset]]$   &   $[[D_2],[D_1],[E_2]]$   & $[[\emptyset],[(0,0)],[(0,0)]]$ &      $D_2+ E_2$ \\ 	


$[[D_1, E_2],[E_2],[\emptyset]]$   &   $[[D_2],[D_1],[E_2]]$   & $[[\emptyset],[(0,0)],[\emptyset]]$ &                $D_2$  \\

$\dots$ & \dots & \dots & \dots  \\

$[[D_1, D_2],[D_2],[\emptyset]]$   &   $[[E_2],[D_2],[\emptyset]]$   & $[[\emptyset],[\emptyset],[(0,0)]]$ &        --- \\ 
$\dots$ & \dots & \dots & \dots  \\\tline
\end{tabular}
\end{table}

When there is no output for $D$ (in the table: ---), it means that there is no divisor corresponding to the particular configuration of generating hyperplanes.

To make the table more readable, we include pictures for the steps marked with a number. 

\vskip .6cm

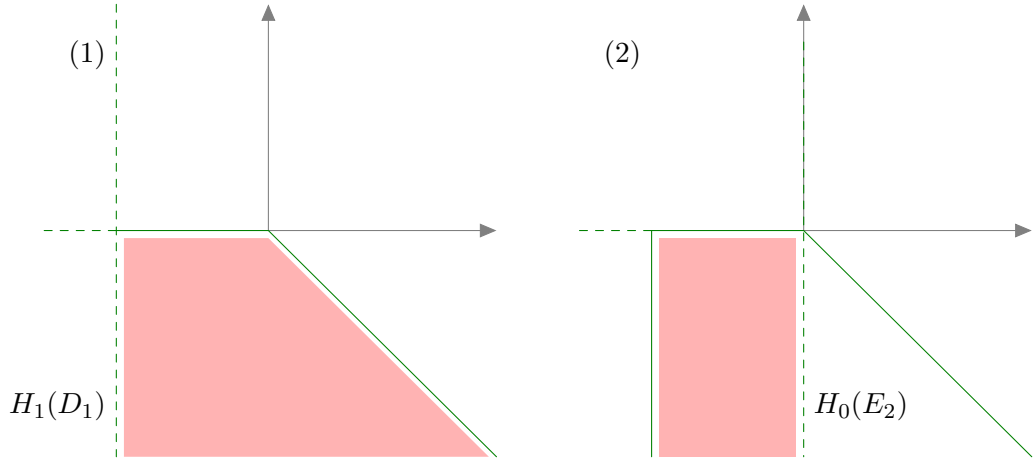
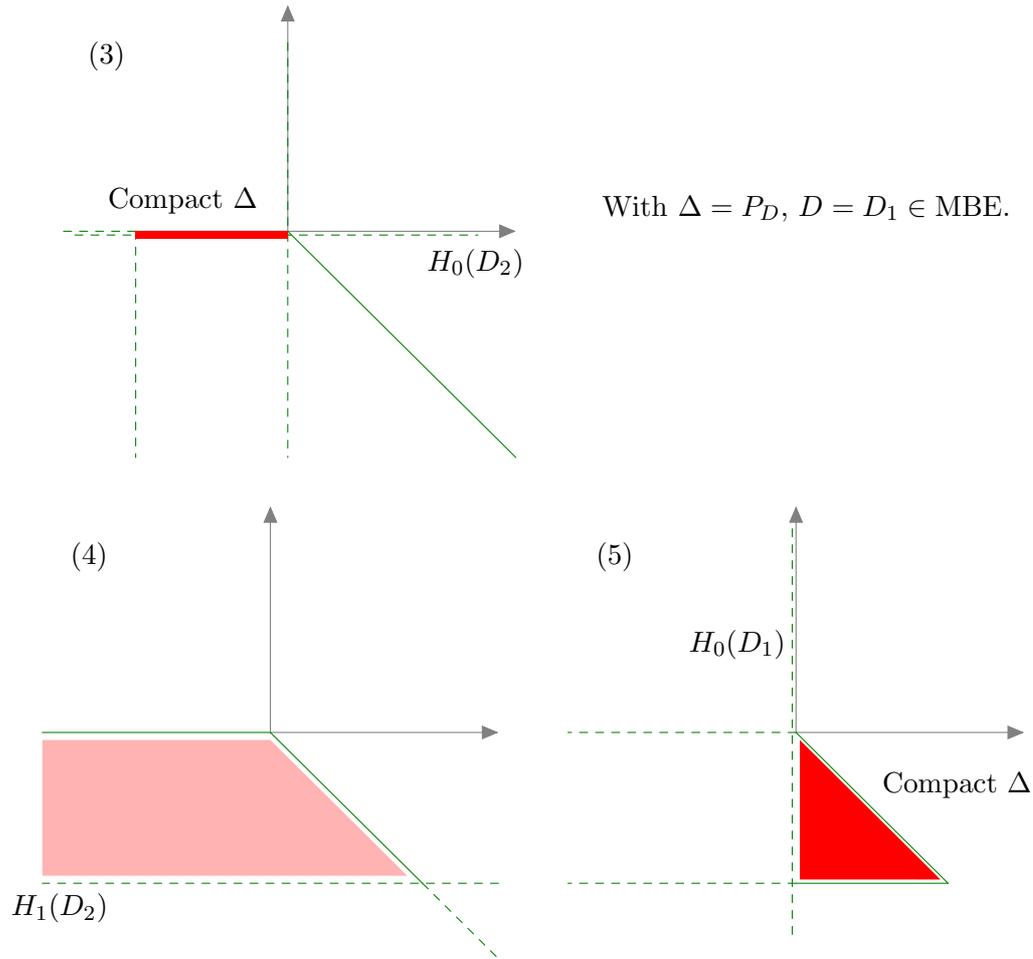
\begin{tikzpicture}[line cap=round,line join=round,>=triangle 45,x=1.0cm,y=1.0cm]
\draw[->,color=gray] (0,0) -- (3,0);
\draw[->,color=gray] (0,0) -- (0,3);
\draw[dashed, color= green!50!black] (-2, 3) -- (-2,-3);
\draw[color=green!50!black] (0,0) -- (-2,0);
\draw[dashed,color=green!50!black] (-2,0) -- (-3,0);
\draw (-2,-2) node [anchor=north east] {$H_1(D_1)$};
\draw[color=green!50!black] (0,0) -- (3,-3);
\fill[color=red,fill=red,fill opacity=0.3] (0,-0.1) -- (2.9,-3) -- (-1.9,-3) -- (-1.9, -0.1) -- cycle;
\draw (-2,2) node [anchor=south east] {(1)};
\end{tikzpicture}
$\qquad$
\begin{tikzpicture}[line cap=round,line join=round,>=triangle 45,x=1.0cm,y=1.0cm]
\draw[->,color=gray] (0,0) -- (3,0);
\draw[->,color=gray] (0,0) -- (0,3);
\draw[color= green!50!black] (-2, 0) -- (-2,-3);
\draw[dashed, color= green!50!black] (0, 2.5) -- (0,-3);
\draw[color=green!50!black] (0,0) -- (-2,0);
\draw[dashed,color=green!50!black] (-2,0) -- (-3,0);
\draw (0,-2) node [anchor=north west] {$H_0(E_2)$};
\draw[color=green!50!black] (0,0) -- (3,-3);
\fill[color=red,fill=red,fill opacity=0.3] (-0.1,-0.1) -- (-0.1,-3) -- (-1.9,-3) -- (-1.9, -0.1) -- cycle;
\draw (-2,2) node [anchor=south east] {(2)};
\end{tikzpicture}

\vskip .6cm

\begin{tikzpicture}[line cap=round,line join=round,>=triangle 45,x=1.0cm,y=1.0cm]
\draw[->,color=gray] (0,0) -- (3,0);
\draw[->,color=gray] (0,0) -- (0,3);
\draw[dashed, color= green!50!black] (-2, 0) -- (-2,-3);
\draw[dashed, color= green!50!black] (0, 2.5) -- (0,-3);
\draw[color=green!50!black] (0,0) -- (-2,0);
\draw[dashed,color=green!50!black] (-2,0) -- (-3,0);
\draw[dashed,color=green!50!black] (-2.8,-0.05) -- (2.5,-0.05);
\draw (1.7,-0.1) node [anchor=north west] {$H_0(D_2)$};
\draw (-2.5,0.1) node [anchor=south west] {Compact $\Delta$};
\draw[color=green!50!black] (0,0) -- (3,-3);
\fill[color=red,fill=red] (0,-0.1) -- (0,0) -- (-2,0) -- (-2, -0.1) -- cycle;
\draw (-2,2) node [anchor=south east] {(3)};
\draw (-3.5,2) node [anchor=south east] {\emph{}};
\draw (4,0) node [anchor=south west] {With $\Delta=P_D$, $D=D_1 \in {\rm MBE}$.};
\end{tikzpicture}

\vskip .6cm

\begin{tikzpicture}[line cap=round,line join=round,>=triangle 45,x=1.0cm,y=1.0cm]
\draw[->,color=gray] (0,0) -- (3,0);
\draw[->,color=gray] (0,0) -- (0,3);
\draw[dashed, color= green!50!black] (3,-2) -- (-3,-2);
\draw[color=green!50!black] (0,0) -- (-3,0);
\draw[dashed,color=green!50!black] (2,-2) -- (3,-3);
\draw (-2,-2) node [anchor=north east] {$H_1(D_2)$};
\draw[color=green!50!black] (0,0) -- (2,-2);
\fill[color=red,fill=red,fill opacity=0.3] (0,-0.1) -- (1.8,-1.9) -- (-3,-1.9) -- (-3, -0.1) -- cycle;
\draw (-2,2) node [anchor=south east] {(4)};
\end{tikzpicture}
\qquad
\begin{tikzpicture}[line cap=round,line join=round,>=triangle 45,x=1.0cm,y=1.0cm]
\draw[->,color=gray] (0,0) -- (3,0);
\draw[->,color=gray] (0,0) -- (0,3);
\draw[color= green!50!black] (2,-2) -- (0,-2);
\draw[dashed, color= green!50!black] (0,-2) -- (-3, -2);
\draw[dashed, color=green!50!black] (0,0) -- (-3,0);
\draw[dashed,color=green!50!black] (-0.05,2.7) -- (-0.05,-2.7);
\draw (0,1.5) node [anchor=north east] {$H_0(D_1)$};
\draw (1,-1) node [anchor=south west] {Compact $\Delta$};
\draw[color=green!50!black] (0,0) -- (2,-2);
\fill[color=red,red] (0.05,-0.1) -- (1.9,-1.95) -- (0.05,-1.95) -- cycle;
\draw (-2,2) node [anchor=south east] {(5)};
\draw (0,-3) node [anchor=south east] {\emph{}};
\end{tikzpicture}

\vskip .5cm

In this case $\Delta=P_D$ with $D=D_2+E_2 \in {\rm MBE}$.

\begin{remark} In the 2-dimensional case there is no need of repeating the algorithm since the nef and the movable cone coincide.
\end{remark}

 \subsection{The blow-up of $\P^3$ in two intersecting lines}
\label{p3-1}
Let us consider the blow-up of $\mathbb{P}^3$ in two intersecting lines. The fan $\Sigma$ contains $6$ rays which are spanned by the points $(1,0,0)$, $(0,1,0)$, $(0,0,1)$, $(-1,-1,-1)$, $(0,-1,-1)$, $(-1,0,-1)$, and correspond to divisor classes $D_1= H-E_1, \; D_2= H-E_2, \; D_3= H,\;  D_4= H-E_1-E_2, \; E_1 \mbox{ and } E_2$ respectively. Let us fix the flag given by $Y_{\bullet}=\{D_1, D_1 \cap D_2, D_1\cap D_2 \cap E_1\}$.
\vskip .5cm
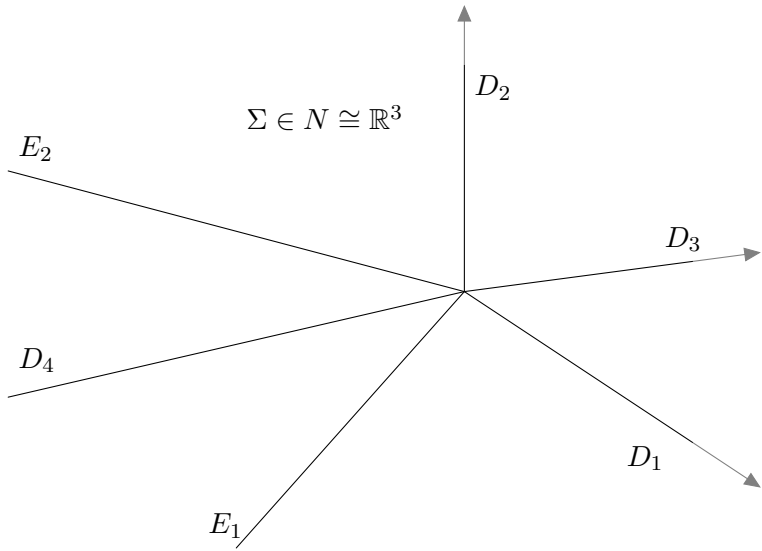
\begin{tikzpicture}[line cap=round,line join=round,>=triangle 45,x=1.0cm,y=1.0cm]
\draw[->,color=gray] (3,0.4) -- (3.9, 0.52);
\draw[->,color=gray] (3,-2) -- (3.9,-2.6);
\draw[->,color=gray] (0,3) -- (0, 3.8);
\draw[color=black] (0,0) -- (3,0.4);
\draw[color=black] (0,0) -- (3,-2);
\draw[color=black] (0,0) -- (0,3);
\draw[color=black] (0,0) -- (-6, -1.4);
\draw[color=black] (0,0) -- (-3, -3.4);
\draw[color=black] (0,0) -- (-6, 1.6);
\draw (2.5,0.4) node [anchor=south west] {$D_3$};
\draw (2,-2.5) node [anchor=south west] {$D_1$};
\draw (0,3) node [anchor=north west] {$D_2$};
\draw (-6,-1.2) node [anchor=south west] {$D_4$};
\draw (-3.5,-3.4) node [anchor=south west] {$E_1$};
\draw (-6,1.6) node [anchor=south west] {$E_2$};
\draw (-3, 2) node [anchor=south west] {$\Sigma \in N \cong \R^3$};
\end{tikzpicture}

\vskip 1.5cm
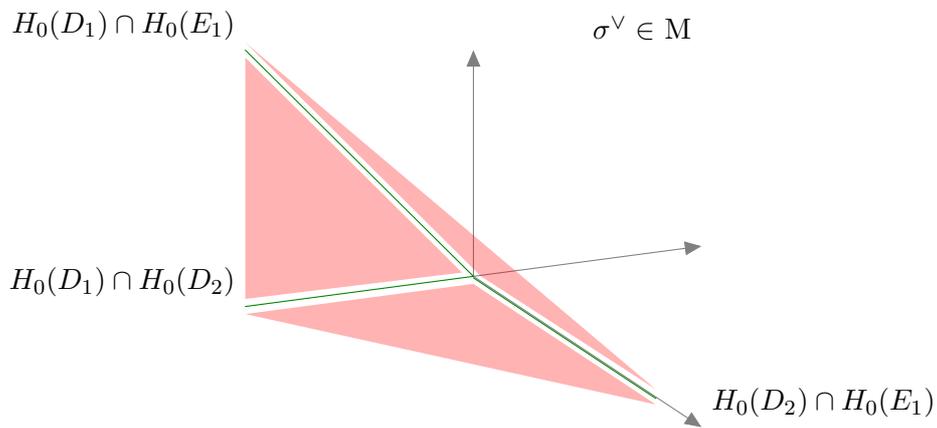
\begin{tikzpicture}[line cap=round,line join=round,>=triangle 45,x=1.0cm,y=1.0cm]
\draw[->,color=gray] (0,0) -- (3, 0.4);
\draw[->,color=gray] (0,0) -- (3,-2);
\draw[->,color=gray] (0,0) -- (0,3);
\draw[color=green!50!black] (0,0) -- (-3, -0.4);
\draw[color=green!50!black] (0,-.02) -- (2.4, -1.62);
\draw[color=green!50!black] (0,0) -- (-3, 3);
\draw (3,-2) node [anchor=south west] {$H_0(D_2)\cap H_0(E_1)$};
\draw (-3,-0.4) node [anchor=south east] {$H_0(D_1)\cap H_0(D_2)$};
\draw (3,3) node [anchor=south east] {$\sigma^{\vee} \in$ M};
\draw (-3,3) node [anchor=south east] {$H_0(D_1)\cap H_0(E_1)$};
\fill[color=red,fill=red,fill opacity=0.3] (-0.15,0.05) -- (-3, 2.9) -- (-3,-0.3) -- cycle;
\fill[color=red,fill=red,fill opacity=0.3] (0,-0.1) -- (-3,-0.5) -- (2.4,-1.7) -- cycle;
\fill[color=red,fill=red,fill opacity=0.3] (0.1,0) -- (-3,3.1) -- (2.4,-1.5) -- cycle;
\end{tikzpicture}

\vskip .5cm
Since this variety admits flips (reversing the order of the blow-ups) and we have chosen a cone that can be flipped, we expect not to find all the generators of the moving cone at the first stage.

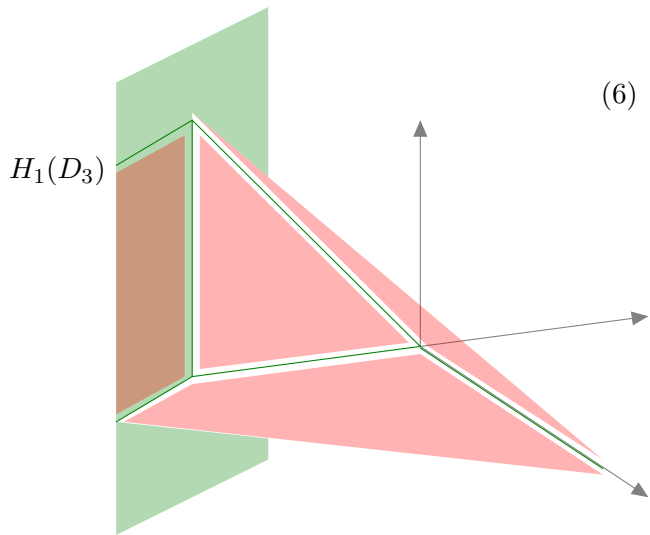
\renewcommand{\arraystretch}{1.3}
\begin{table}[ht]\centering\scriptsize
      \begin{tabular}[t]{c|c|c|c}\tline
         $R $ &  $M$   &$P$ &       D \\\tline
 	$[[D_3,D_4, E_2],[D_3, D_4,E_2],[D_3, D_4,E_2]]$   &   $[[\emptyset],[\emptyset],[\emptyset]]$   & $[\emptyset]$ &        \\ 	

	
	$[[D_4, E_2],[D_4,E_2],[D_4,E_2]]$   &   $[[D_3],[\emptyset],[\emptyset]]$   & $[[\emptyset], [\vec0, (0, 0, -1), (0,1,-1)], [\emptyset]]$ &              $ -$ (6)\\ 	
	
$\dots$ & \dots & \dots & \dots  \\

	$[[D_4, E_2],[E_2],[\emptyset]]$   &   $[[D_3],[D_4],[E_2]]$   & $[[\emptyset],[\vec0],[\vec0]]$ &        $D_3$  \\
	$[[D_4, E_2],[E_2],[\emptyset]]$   &   $[[D_3],[D_4],[E_2]]$   & $[[\emptyset],[\vec0],[\emptyset]]$ &        $D_3$  \\

$\dots$ & \dots & \dots & \dots  \\

	$[[D_4, E_2],[E_2],[\emptyset]]$   &   $[[D_3],[D_4],[E_2]]$   & $[[\emptyset],[\emptyset],[\vec0]]$ &        $D_3$  \\
	$[[D_4, E_2],[E_2],[\emptyset]]$   &   $[[D_3],[D_4],[E_2]]$   & $[[\emptyset],[\emptyset],[\emptyset]]$ &        $D_3$  \\

$\dots$ & \dots & \dots & \dots  \\

	$[[D_4, E_2],[D_4],[\emptyset]]$   &   $[[D_3],[E_2],[D_4]]$   & $[[\emptyset],[\vec0],[\vec0]]$ &       $D_4 + D_3$  \\
	
	$[[D_4, E_2],[D_4],[\emptyset]]$   &   $[[D_3],[E_2],[D_4]]$   & $[[\emptyset],[\vec0],[\emptyset]]$ &       $D_3$  \\

$\dots$ & \dots & \dots & \dots  \\

	$[[D_4, E_2],[D_4],[\emptyset]]$   &   $[[D_3],[E_2],[D_4]]$   & $[[\emptyset],[\emptyset],[\vec0]]$ &       $D_4 + D_3$  \\
	
	$[[D_4, E_2],[D_4],[\emptyset]]$   &   $[[D_3],[E_2],[D_4]]$   & $[[\emptyset],[\emptyset],[\emptyset]]$ &       $D_3$  \\





$\dots$ & \dots & \dots & \dots  \\
$[[D_3, E_2],[E_2],[\emptyset]]$   &   $[[D_4],[D_3],[E_2]]$   & $[[\emptyset],[\vec0],[\vec0]]$ &      $D_4+ E_2$ \\ 	


$[[D_3, E_2],[E_2],[\emptyset]]$   &   $[[D_4],[D_3],[E_2]]$   & $[[\emptyset],[\vec0],[\emptyset]]$ &                ---  \\

$\dots$ & \dots & \dots & \dots  \\\tline
\end{tabular}
\end{table}

\begin{tikzpicture}[line cap=round,line join=round,>=triangle 45,x=1.0cm,y=1.0cm]
\draw[->,color=gray] (0,0) -- (3, 0.4);
\draw[->,color=gray] (0,0) -- (3,-2);
\draw[->,color=gray] (0,0) -- (0,3);
\draw[color=green!50!black] (0,0) -- (-3, -0.4);
\draw[color=green!50!black] (0,-.02) -- (2.4, -1.62);
\draw[color=green!50!black] (0,0) -- (-3, 3);
\draw (-4,2) node [anchor=south east] {$H_1(D_3)$};
\fill[color=green!50!black,fill=green!50!black,fill opacity=0.3] (-4, 3.5) -- (-4, -1)-- (-3, -0.4 ) -- (-3, 3.1) -- (-2, 2.25) -- (-2, 4.5)-- cycle;
\fill[color=green!50!black,fill=green!50!black,fill opacity=0.3] (-4, -1)--(-4, -2.5)--(-2, -1.5 ) -- (-2, -1.22)-- cycle;
\fill[color=red,fill=red,fill opacity=0.3] (-0.15,0.05) -- (-2.9, 2.8) -- (-2.9,-0.3) -- cycle;
\fill[color=red,fill=red,fill opacity=0.3] (0,-0.1) -- (-3,-0.5) -- (-3.9, -1) -- (2.4,-1.7) -- cycle;
\fill[color=red,fill=red,fill opacity=0.3] (0.1,0) -- (-3,3.1) -- (2.4,-1.5) -- cycle;
\fill[color=red,fill=red,fill opacity=0.3] (-3.1, 2.8) -- (-3.1,-0.4)-- (-4,-0.9) -- (-4, 2.3)  -- cycle; 
\draw[color=green!50!black] (-3,-0.4) -- (-3, 3);
\draw[color=green!50!black] (-3,-0.4) -- (-4, -1);
\draw[color=green!50!black] (-4,2.4) -- (-3, 3);
\draw (3,3) node [anchor=south east] {(6)};
\end{tikzpicture}

\begin{remark}
Since we have chosen a flag based on a cone that can be flipped, there are two special features appearing:
\begin{itemize}
\item There are extremal rays of the movable cone missing. In fact not every movable divisor will have $\chi^0$ as local section, whereas this is true for nef divisors. This will be remedied by looking at another $T$-invariant flag and completing with the new Minkowski basis elements found.
\item The divisor $D_4+D_3 = 2H - E_1 -E_2$ still corresponds to an irreducible polytope, since it is only decomposable by movable divisors $2H - E_1 -E_2 = (H-E_1)+(H-E_2)$ and additivity of polytopes is only respected by nef divisors.
\end{itemize}
\end{remark}

\subsection{Flip invariant flag}
\label{p3-2}
In the setting of the previous example, we see that choosing a flag not contained in the base locus of 
any effective divisor we obtain all the generators just in one step.

\vskip .2cm
Choosing the flag $Y_{\bullet}=\{D_1, D_1 \cap D_2, D_1\cap D_2 \cap D_3\}$ we obtain all basis elements.

\renewcommand{\arraystretch}{1.3}
\begin{table}[ht]\centering\footnotesize
      \begin{tabular}[t]{c|c|c|c}\tline
         $R $ &  $M$   &$P$ &       D \\\tline

	
	

$\dots$ & \dots & \dots & \dots  \\

	$[[E_1, E_2],[E_2],[\emptyset]]$   &   $[[D_4],[E_1],[E_2]]$   & $[[\emptyset],[\vec0],[\vec0]]$ &        $D_4+E_1+E_2$  \\

	$[[E_1, E_2],[E_2],[\emptyset]]$   &   $[[D_4],[E_1],[E_2]]$   & $[[\emptyset],[\vec0],[\emptyset]]$ &        $D_4+E_1$  \\

$\dots$ & \dots & \dots & \dots  \\
	$[[E_1, E_2],[E_2],[\emptyset]]$   &   $[[D_4],[E_1],[E_2]]$   & $[[\emptyset],[\emptyset],[\vec0]]$ &        $D_4+E_2$  \\
	$[[E_1, E_2],[E_2],[\emptyset]]$   &   $[[D_4],[E_1],[E_2]]$   & $[[\emptyset],[\emptyset],[\emptyset]]$ &       ---  \\

$\dots$ & \dots & \dots & \dots \\

	
	


	$[[D_4, E_2],[D_4],\emptyset]]$   &   $[[E_1],[E_2],[D_4]]$   & $[[\emptyset],[\vec0], [\vec0]]$ &              $2D_4+E_1+E_2$\\ 

$\dots$ & \dots & \dots & \dots  \\\tline
\end{tabular}
\end{table}



\section{Finding the decomposition}\label{decompose}

Let us now turn to the problem of finding the Minkowski decomposition of a given $T$-invariant big divisor $D$ once the Minkowski basis for $X$ is determined by the algorithm. 

Two different equivalent approaches are possible: we can either determine the cone of the secondary fan in which the movable part $M_D$ of $D$ lies and then write $M_D$ as a non-negative linear combination of its extremal rays, or we can invest knowledge of the Okounkov body $\Delta_{Y_\bullet}(D)$ with respect to some $T$-invariant flag and go backwards by finding Minkowski summands of this polytope. Note that the first option can be quite challenging since it presumes knowledge of the whole structure of the secondary fan as well as of the decomposition of $D$ into fixed and movable part. On the other hand the input for the second option comes naturally, we just have to describe how to find the Minkowski summands of a polytope, but this is given by Proposition \ref{ref} as follows.

Given the divisor $D$, we consider its polytope $P_D$ and the corresponding dual fan $\Sigma_D$. Now, a Minkowski basis element $B$ is a summand in the Minkowski decomposition of $D$ if and only if its fan $\Sigma_B$ is refined by $\Sigma_D$. This is a straightforward check. Note also that the fans $\Sigma_D$ and $\Sigma_B$ only dependent on the linear equivalence classes. 
We thus obtain all classes of Minkowski basis elements with representative $B_i$ having positive coefficient $a_i$ in some Minkowski decomposition $M_D=\sum a_jB_j$. Note that from this data we can read off important information that in general can be hard to obtain: in particular for a moving divisor $D$ the above procedure immediately tells us in which cone of the secondary fan $D$ lies. 

\newpage
\begin{example} 
Let us consider $X$ the blow-up of $\P^2$ in two points. We denote the exceptional divisors with $E_1,E_2$ and the class of the pullback of a general line by $H$. The corresponding fan is spanned by the following rays in $\R^2$.
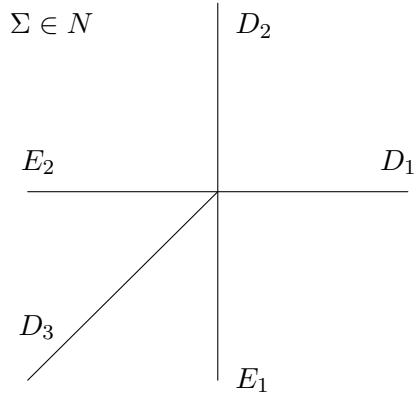
\emph{}
\begin{center}
\begin{tikzpicture}[line cap=round,line join=round,>=triangle 45,x=1.0cm,y=1.0cm]
\draw[color=black] (-2.5,0) -- (2.5,0);
\draw (0.1,2.5) node [anchor=north west] {$D_2$};
\draw (2,0.1) node[anchor=south west] {$D_1$};
\draw (0.1, -2.5) node[anchor=west] {$E_1$};
\draw (-2,-1.5) node[anchor=north east] {$D_3$};
\draw (-2,0.1) node[anchor=south east] {$E_2$};
\draw[color=black] (0,-2.5) -- (0,2.5);
\draw[color=black] (0,0) -- (-2.5,-2.5);
\draw (-1.5,2.5) node[anchor=north east] {$\Sigma \in N$};
\end{tikzpicture}
\end{center}
It is easy to see that the nef cone $\Nef(X)$ is spanned by the classes $M_1 := H$, $M_2:=H-E_1$ and $M_3 := H-E_2$. By Theorem \ref{th:basis} these classes form a Minkowski basis with respect to torus-invariant flags. The normal fans of their polytopes are depicted below.\\

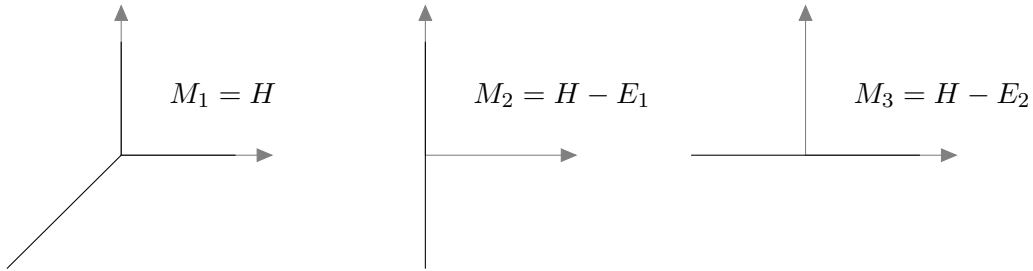
\begin{tikzpicture}[line cap=round,line join=round,>=triangle 45,x=1.0cm,y=1.0cm]
\draw[->,color=gray] (0,0) -- (2,0);
\draw[->,color=gray] (0,0) -- (0,2);

\draw (0,0) -- (1.5,0);
\draw (0,0) -- (0,1.5);
\draw (0,0) -- (-1.5,-1.5);

\draw (.5,.5) node[anchor=south west] {$M_1=H$};

\draw[->,color=gray] (4,0) -- (6,0);
\draw[->,color=gray] (4, 0) -- (4,2);
\draw (4, -1.5) -- (4,1.5);

\draw (4.5,.5) node[anchor=south west] {$M_2=H-E_1$};
\draw[->,color=gray] (9,0) -- (11,0);
\draw[->,color=gray] (9, 0) -- (9,2);
\draw (7.5,0) -- (10.5,0);
\draw (9.5,0.5) node[anchor=south west] {$M_3=H-E_2$};

\end{tikzpicture}\\

Let us now decompose the $T$-invariant divisor $D:= D_1+D_2+E_1$. It has the following polytope $P_D$ and corresponding normal fan $\Sigma_{P_D}$.\\

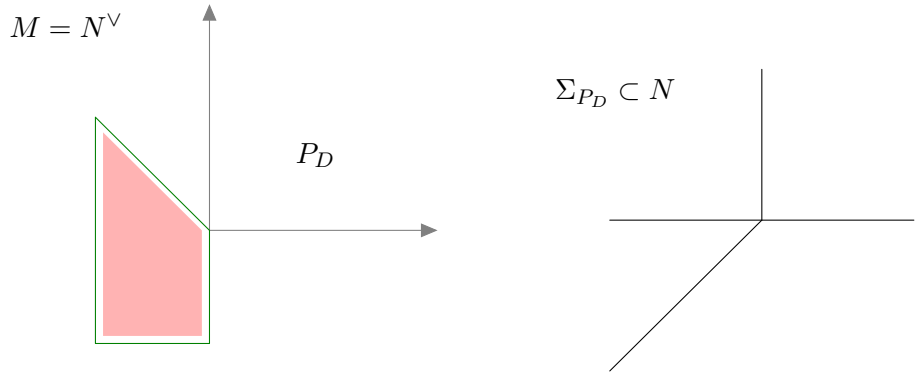
\begin{tikzpicture}[line cap=round,line join=round,>=triangle 45,x=1.0cm,y=1.0cm]
\draw[->,color=gray] (0,0) -- (3,0);
\draw[->,color=gray] (0,0) -- (0,3);
\draw[color=green!50!black] (-1.5,1.5) -- (-1.5,-1.5);
\draw[color=green!50!black] (-1.5,-1.5) -- (0,-1.5);
\draw[color=green!50!black] (0,-1.5) -- (0,0);
\draw[color=green!50!black] (0,0) -- (-1.5,1.5);
\fill[color=red,fill=red,fill opacity=0.3] (-1.4,1.3) -- (-1.4,-1.4) -- (-.1,-1.4) -- (-.1, 0 ) -- cycle;
\draw (1, 1) node[anchor=west] {$P_D$};
\draw (1, -2.5) node[anchor=west] {\emph{}};
\draw (-1,3) node[anchor=north east] {$M=N^{\vee}$};
\draw (-1,-2) node[anchor=north east] {\emph{}};
\end{tikzpicture}
$\quad \quad \quad$ 
\begin{tikzpicture}[line cap=round,line join=round,>=triangle 45,x=1.0cm,y=1.0cm]
\draw[color=black] (-2,0) -- (2,0);
\draw[color=black] (0,0) -- (0,2);
\draw[color=black] (0,0) -- (-2,-2);
\draw (-1,2) node[anchor=north east] {$\Sigma_{P_D} \subset N$};
\draw (-1,-2.5) node[anchor=north east] {\emph{}};
\end{tikzpicture}

Now, $\Sigma_{P_D}$ is the refinement of the fans $\Sigma_{M_1}$ and $\Sigma_{M_3}$ but not of $\Sigma_{M_2}$. Thus $D$ decomposes as a positive linear combination of $M_1$ and $M_3$. We easily obtain the decomposition $D \sim M_1+M_3 = H + (H-E_2)$. 

\end{example}



\vskip .2cm

 Piotr Pokora, 
Institute of Mathematics, 
Pedagogical University of Cracow, 
Podchor\c{a}\.zych 2, 
30-084 Krak\'ow, Poland.

 \nopagebreak
   \noindent \textit{E-mail address:} \texttt{piotrpkr@gmail.com}
 
  \vskip .2cm
 David Schmitz,
   Fach\-be\-reich Ma\-the\-ma\-tik und In\-for\-ma\-tik,
   Philipps-Uni\-ver\-si\-t\"at Mar\-burg,
   Hans-Meer\-wein-Stra{\ss}e,
   D-35032~Mar\-burg, Germany.

 \nopagebreak
\noindent   \textit{E-mail address:} \texttt{schmitzd@mathematik.uni-marburg.de}

  \vskip .2cm
   Stefano Urbinati,
   Universit\`a degli Studi di Padova,
   Dipartimento di Matematica,
   ROOM 608, Via Trieste 63,
   35121 Padova, Italy.

   \nopagebreak
\noindent   \textit{E-mail address:} \texttt{urbinati.st@gmail.com}


\end{document}